\documentclass[12pt]{article}
\usepackage[utf8]{inputenc}
\usepackage[left=1in,right=1in, top=1.2in,bottom=1.2in,includehead]{geometry}
\usepackage{amssymb,amsthm,amsmath,pgf}
\usepackage{mathrsfs}
\usepackage{placeins}
\usepackage{wrapfig}
\usepackage{indentfirst}
\usepackage{fancyhdr}
\usepackage{url}

\usepackage[hidelinks]{hyperref}

\newtheorem{theorem}{Theorem}[section]
\newtheorem{lemma}{Lemma}[section]

\newtheorem{corollary}{Corollary}[section]

\newtheorem{definition}{Definition}

\newcommand{\newcase}[1]{\vspace{0.7em}\noindent\textbf{Case #1}\vspace{0.5em}\par}
\newcommand{\cbeginproof}[0]{\par\noindent\textit{Proof.} }
\newcommand{\cendproof}[0]{ \qed\par\vspace{1em}}

\begin{document}

\title{Fault-Tolerant Detection Systems on the King's Grid}
\author{
    \small Devin C. Jean\\
    \small Vanderbilt University\\
    \small \texttt{devin.c.jean@vanderbilt.edu}
    \and
    \small Suk J. Seo\\
    \small Middle Tennessee State University\\
    \small \texttt{Suk.Seo@mtsu.edu}
}
\date{}

\maketitle
\thispagestyle{empty}

\begin{abstract}
A detection system, modeled in a graph, uses ``detectors" on a subset of vertices to uniquely identify an ``intruder" at any vertex. We consider two types of detection systems: open-locating-dominating (OLD) sets and identifying codes (ICs).
An OLD set gives each vertex a unique, non-empty open neighborhood of detectors, while an IC provides a unique, non-empty closed neighborhood of detectors.
We explore their fault-tolerant variants: redundant OLD (RED:OLD) sets and redundant ICs (RED:ICs), which ensure that removing/disabling at most one detector guarantees the properties of OLD sets and ICs, respectively.
This paper focuses on constructing optimal RED:OLD sets and RED:ICs on the infinite king's grid, and presents the proof for the bounds on their minimum densities;   $\left[\frac{3}{10}, \frac{1}{3}\right]$ for RED:OLD sets and $\left[\frac{3}{11}, \frac{1}{3}\right]$  for RED:ICs.
\end{abstract}

\vspace{1em}


\section{Introduction}

Let $G = (V(G), E(G))$ be an (undirected) graph, with vertices $V(G)$ and edges $E(G)$, modelling a system or facility with detectors to recognize a possible problem, traditionally referred to as an ``intruder".
For example, the vertices of the graph can represent sections of a shopping mall, the intruder could be a shoplifter, and the detectors can be video surveillance equipment or motion, magnetic, or RFID sensors.
The goal here is to identify the exact location/vertex of the intruder by placing the minimum number of detectors in the facility/graph. 

A \emph{distinguishing set} $\Re \subseteq \mathscr{P}(V(G))$ is a collection of subsets of vertices which covers $V(G)$ and provides the ability to separate any vertex from the rest of the vertices, and therefore it can be used to pinpoint the location of the intruder.
A more practical method, however, is to consider a subset of vertices $S \subseteq V(G)$ called a ``detection system", where each vertex in $S$ is installed with a specific type of detector or sensor for locating the intruder.  

\begin{definition}\label{def:det-region}
For detector set $S \subseteq V(G)$, the \emph{detection region} of $v \in S$, denoted $R_{det}(v) \subseteq V(G)$, is the region in which $v$ can detect an intruder.
\end{definition}

In a detection system, the detector installed at vertex $v$ is capable of detecting an intruder if and only if it is within $R_{det}(v)$. Note that $R_{det}(v)$ has also been referred to as the ``watching zone" of $v$ in other papers, but with different notation \cite{watchsys-1,watchsys-2}. 
Using an appropriate selection of detectors $S \subseteq V(G)$ for a graph, we can generate $\{R_{det}(v) : v \in S\}$, which forms a distinguishing set.
Ideally, we would like to use the smallest number of detectors possible while managing to generate a distinguishing set for $V(G)$ that can be used to pinpoint the location of the intruder.

\begin{definition}\cite{oldking}
Let $G$ be a graph and $v \in V(G)$.
The \emph{open neighborhood} of $v$, denoted $N(v)$, is the set of all vertices adjacent to $v$, $\{w \in V(G) : vw \in E(G)\}$.
\end{definition}

\begin{definition}
Let $G$ be a graph and $v \in V(G)$.
The \emph{closed neighborhood} of $v$, denoted $N[v]$, is the set of all vertices adjacent to $v$ as well as $v$ itself, $N(v) \cup \{v\}$.
\end{definition}

Many types of detection systems with various properties have been explored throughout the years.
Over 470 papers have been published on various types of detector-based sets, distinguishing sets, and other related concepts \cite{dombib}.
One such system is the \emph{Locating-Dominating (LD) set}, where each detector identifies an intruder in their closed neighborhood and also has the ability to distinguish the vertex itself from its neighbors \cite{ftld,dom-loc-acyclic}.
Of particular interest in this paper are \emph{open-locating-dominating (OLD)} sets \cite{oldking}, which are detection systems using $R_{det}(v) = N(v)$, and \emph{identifying codes (IC)} \cite{karpovsky}, which are detection systems using $R_{det}(v) = N[v]$.
Note that both OLD sets and ICs generate distinguishing sets which cover the entire graph using detectors that detect intruders at distance at most 1.
The only difference is that OLD detectors cannot detect an intruder at their own vertex while IC detectors must do so.
For examples of minimum a OLD set and IC on the same graph, see Figure~\ref{fig:ex-finite} (a) and (c), respectively.

An OLD set is a form of \emph{open-dominating set}, which simply requires that $\cup_{v \in S}{N(v)} = V(G)$ where $S \subseteq V(G)$ is the open-dominating set \cite{oldking}.
However, open-domination alone is insufficient to generate a distinguishing set, which is why OLD sets add additional uniqueness requirements in order to distinguish vertices.
Similarly, an IC is a form of \emph{dominating set}, which simply requires that $\cup_{v \in S}{N[v]} = V(G)$ where $S \subseteq V(G)$ is the dominating set.

In the next section, we introduce the general concept of a \emph{detection system}---which is similar to a ``watching system" as used in other papers \cite{watchsys-1,watchsys-2}---and fault-tolerant variants known as a \emph{$k$-redundant detection systems}.
In Section~\ref{sec:lower-bound-common} we introduce the concept of a share argument, which will be used for later proofs.
Finally, in Sections \ref{sec:red-old} and \ref{sec:red-ic} we explore upper and lower bounds for optimal 1-redundant variants of OLD sets and ICs known as RED:OLD sets and RED:ICs, respectively.

\section{Detection Systems}

In a detector-based distinguishing set we use what is known as the ``locating code" of a vertex to identify the intruder based on which sensors have detected the intruder in their detection region.

\begin{definition}\label{def:loc-code}
Given a set of detectors $S \subseteq V(G)$ and a vertex $v \in V(G)$, the \emph{locating code} of $v$, denoted $\mathscr{L}_S(v)$, is the set of detectors which dominate (cover) $v$. $\mathscr{L}_S(v) \equiv \{w \in S : v \in R_{det}(w)\}$.
\end{definition}

The detection region function $R_{det}$ is only required to be defined over the detector set $S \subseteq V(G)$.
However, if it defined for every vertex---that is, a hypothetical detection region for a detector being placed at any given vertex---then it is said to be \emph{uniform}.

\begin{definition}
$R_{det}$ is said to be \emph{uniform} if it is defined for every vertex in $V(G)$.
\end{definition}

\begin{definition}
$R_{det}$ is said to be \emph{symmetric} if $u \in R_{det}(v) \implies v \in R_{det}(u)$.
\end{definition}

\begin{theorem}\label{theo:loc-code}
Let $S \subseteq V(G)$ be a set of detectors.
If $R_{det}$ is uniform and symmetric then for any $v \in V(G)$, $\mathscr{L}_S(v) = R_{det}(v) \cap S$.
\end{theorem}
\begin{proof}
$\mathscr{L}_S(v) \equiv \{w \in S : v \in R_{det}(w)\} = \{w \in S : w \in R_{det}(v)\} = R_{det}(v) \cap S$
\end{proof}

Theorem~\ref{theo:loc-code} is not strictly necessary to compute locating codes, but when applicable it can be more convenient than Definition~\ref{def:loc-code}.
Notably, OLD sets use $R_{det}(v) = N(v)$ and ICs use $R_{det}(v) = N[v]$ for every detector regardless of its position; thus they are uniform.
And for undirected graphs---which this paper assumes---$N(v)$ and $N[v]$ are also symmetric.
Thus, for OLD sets and ICs, $\mathscr{L}_S(v) = R_{det}(v) \cap S$.

Depending on the problem, we might need the detection system to be fully functional even when one of the detectors goes offline or is broken.
For this purpose, we could instead use a fault-tolerant variant of a detection system known as a $k$-redundant detection system, which is a detection system for which removing at most $k$ detectors is still a detection system.
In this paper, we will explore \emph{redundant open-locating-dominating (RED:OLD)} sets and \emph{redundant identifying codes (RED:ICs)}, which are 1-redundant detection systems \cite{our13, ft-old-cubic}.

\begin{definition}\label{def:det-sys}
A \emph{detection system} is a set of vertices $S \subseteq V(G)$ such that $\forall v \in V(G)$, $\mathscr{L}_S(v) \neq \varnothing$ and $\forall u,v \in V(G)$ with $u \neq v$,\; $\mathscr{L}_S(u) \neq \mathscr{L}_S(v)$.
\end{definition}

\begin{definition}\label{def:red-det-sys}
A \emph{k-redundant detection system} is a detection system $S \subseteq V(G)$ such that for any $D \subseteq S$ with $|D| \le k$, $S-D$ is a detection system.
\end{definition}

By Definition~\ref{def:red-det-sys}, a detection system is technically a 0-redundant detection system.
Next, we present characterizations that will be useful in constructing detection systems.

\begin{definition}\label{def:dom}
A vertex $v \in V(G)$ is said to be \emph{k-dominated} by a set $S \subseteq V(G)$ if it is dominated (covered) by the detection regions of exactly $k$ vertices in $S$.
That is, $|\mathscr{L}_S(v)| = k$.
\end{definition}

\begin{definition}\cite{ourtri}\label{def:disty}
A set $S \subseteq V(G)$ is said to be \emph{k-distinguishing} if for any two distinct vertices $u,v \in V(G)$, $|\mathscr{L}_S(u) \triangle \mathscr{L}_S(v)| \ge k$.
\end{definition}

In this paper, when quantifiers such as ``at least" are used with $k$-dominated, they apply to the value of $k$.
For instance, ``at least 1-dominated" means $k$-dominated for some $k \ge 1$.

\begin{definition}\label{def:dom-count}
Let $S \subseteq V(G)$ be a set of detectors.
The \emph{domination count} of a vertex $v \in V(G)$, denoted $dom_S(v)$, is $k$ such that $v$ is $k$-dominated by $S$.
\end{definition}

The notation $dom_S(v)$ introduced by Definition~\ref{def:dom-count} is technically equivalent to $|\mathscr{L}_S(v)|$ but is used in this paper for semantic clarity when only the count is needed.

\begin{theorem}\label{theo:det-sys}
A set $S \subseteq V(G)$ is a detection system if and only if all vertices are at least 1-dominated and all pairs are 1-distinguished.
\end{theorem}
\begin{proof}\label{theo:det-sys-proof}
Consider the properties given in Definition~\ref{def:det-sys} and apply the following equivalences:
\begin{equation*}
\mathscr{L}_S(v) \neq \varnothing \Longleftrightarrow |\mathscr{L}_S(v)| \ge 1
\end{equation*}
\begin{equation*}
\mathscr{L}_S(v) \neq \mathscr{L}_S(u) \Longleftrightarrow \mathscr{L}_S(v) \triangle \mathscr{L}_S(u) \neq \varnothing \Longleftrightarrow |\mathscr{L}_S(v) \triangle \mathscr{L}_S(u)| \ge 1
\end{equation*}
\end{proof}

\begin{theorem}\cite{ftsets}\label{theo:red-det-sys}
A set $S \subseteq V(G)$ is a $k$-redundant detection system ($k \ge 0$) if and only if all vertices are at least $(k+1)$-dominated and all pairs are $(k+1)$-distinguished.
\end{theorem}
\begin{proof}
Let $S \subseteq V(G)$ be a $k$-redundant detection system and $v \in V(G)$.
Suppose $|\mathscr{L}_S(v)| \le k$.
Consider the new set $S' = S-\mathscr{L}_S(v)$.
Then $\mathscr{L}_{S'}(v) = \varnothing$, implying $S'$ is not a detection system, a contradiction.
Thus, $|\mathscr{L}_S(v)| \ge k + 1$ and we have that $S$ at least $(k+1)$-dominates all vertices.
Let $u \in V(G)$ with $u \neq v$ and let $D = \mathscr{L}_S(u) \triangle \mathscr{L}_S(v)$.
Suppose that $|D| \le k$.
Consider the new set $S' = S-D$.
We see that $\mathscr{L}_{S'}(u) \triangle \mathscr{L}_{S'}(v) = \varnothing$, so $S'$ is not a detection system, a contradiction.
Thus, $|D| \ge k + 1$ and we have that $S$ is $(k+1)$-distinguishing.

For the converse, suppose $S \subseteq V(G)$ is $(k+1)$-distinguishing and it at least $(k+1)$-dominates all vertices.
These properties satisfy the requirements of Theorem~\ref{theo:det-sys}, so we know that $S$ is a detection system.
By hypothesis, $S$ is $(k+1)$-distinguishing, so $\forall u,v \in V(G)$ with $u \neq v$, $|\mathscr{L}_S(u) \triangle \mathscr{L}_S(v)| \ge k + 1$.
Additionally, $S$ is at least $(k+1)$-dominating, so $\forall v \in V(G),\; |\mathscr{L}_S(v)| \ge k + 1$.
Let $D$ be an arbitrary subset of $S$ with $|D| \le k$.
Consider the new set $S' = S-D$.
The removal of $|D|$ elements from $S$ will remove at most $|D|$ elements from any locating code; thus, $\forall v \in V(G),\; |\mathscr{L}_{S'}(v)| \ge 1$, implying $S'$ at least 1-dominates all vertices.
Additionally, the removal of $|D|$ elements from $S$ will remove at most $|D|$ differences in any pair of locating codes, thus $\forall u,v \in V(G)$ with $u \neq v$, $|\mathscr{L}_{S'}(u) \triangle \mathscr{L}_{S'}(v)| \ge 1$, so $S'$ is 1-distinguishing.
Therefore, by Theorem~\ref{theo:det-sys}, $S'$ is a detection system, implying that $S$ is a $k$-redundant detection system.
\end{proof}

Note that Theorem~\ref{theo:red-det-sys} has been introduced and proven for the specific case where $k=1$ \cite{ftsets}, but we have provided a generalized form here.

\begin{theorem}\label{theo:det-sys-superset}
Suppose $S \subseteq V(G)$ at least $p$-dominates all vertices and is $q$-distinguishing.
Then any superset $S' \supseteq S$ also at least $p$-dominates all vertices and is $q$-distinguishing.
\end{theorem}
\begin{proof}
Let $T \subseteq V(G)-S$ with $T \neq \varnothing$.
Let $S' = S \cup T$.
Because $S$ at least $p$-dominates all vertices, $\forall v \in V(G),\; |\mathscr{L}_S(v)| \ge p$.
By Definition~\ref{def:loc-code}, $\mathscr{L}_S(v) \equiv \{x \in S : v \in R_{det}(x)\}$.
Therefore, the inclusion of vertices $T$ as new detectors results in $\mathscr{L}_S(v) \subseteq \mathscr{L}_{S'}(v)$, so $|\mathscr{L}_{S'}(v)| \ge p$ and we have that $S'$ at least $p$-dominates all vertices.
Additionally, because $S$ is $q$-distinguishing, for any two distinct vertices $u,v \in V(G),\; |\mathscr{L}_S(u) \triangle \mathscr{L}_S(v)| \ge q$.
By the previous finding, we see that $\mathscr{L}_S(u) \triangle \mathscr{L}_S(v) \subseteq \mathscr{L}_{S'}(u) \triangle \mathscr{L}_{S'}(v)$, so $|\mathscr{L}_{S'}(u) \triangle \mathscr{L}_{S'}(v)| \ge q$ and we have that $S'$ is $q$-distinguishing.
\end{proof}

\begin{corollary}
Any superset of a $k$-redundant detection system is also a $k$-redundant detection system.
\end{corollary}

Disregarding the specific choice of $S'$ used in Theorem~\ref{theo:det-sys-superset}, the theorem gives us an important fact: if there exists a non-trivial detection system or $k$-redundant detection system on some graph then there also exists a larger one.
Therefore, what we are interested in is finding the smallest possible set on a given graph.

For finite graphs, we use the notations OLD($G$), RED:OLD($G$), IC($G$), and RED:IC($G$) to denote the cardinality of the smallest possible such sets on graph $G$, respectively.
However, for infinite graphs, we measure via the \emph{density} of the subset, which is defined as the ratio of the size of the subset to the size of the whole set \cite{ourtri,oldking,ftsets}.
We use the notations OLD\%($G$), RED:OLD\%($G$), IC\%($G$), and RED:IC\%($G$) to denote the lowest density of any possible such set on $G$ \cite{ourtri,oldking,ftsets}.
Note that density is also defined for finite graphs.

\begin{figure}[ht]
    \centering
    \begin{tabular}{c@{\hskip 2em}c@{\hskip 2em}c@{\hskip 2em}c}
        \includegraphics[width=0.19\textwidth]{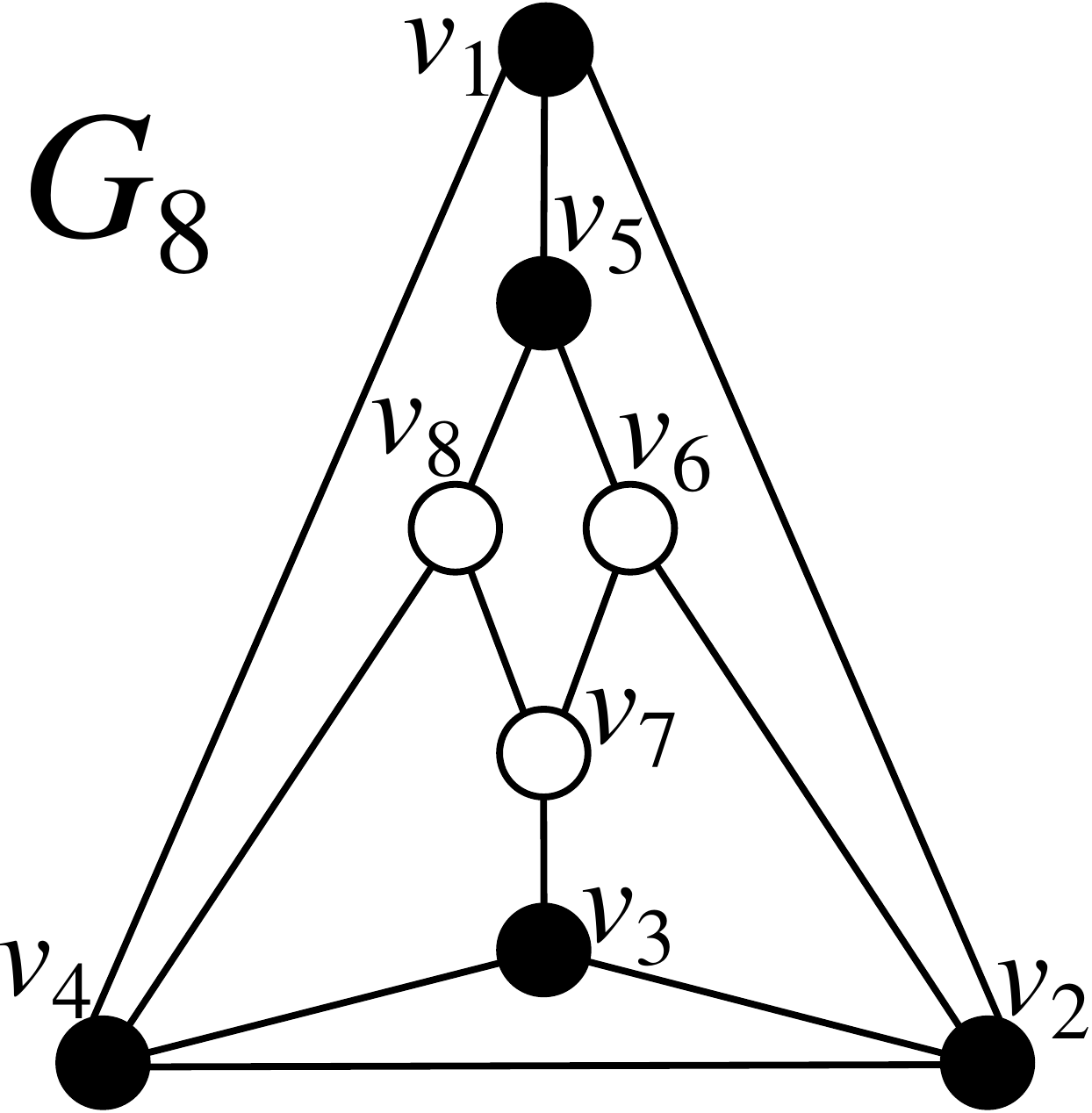} & \includegraphics[width=0.19\textwidth]{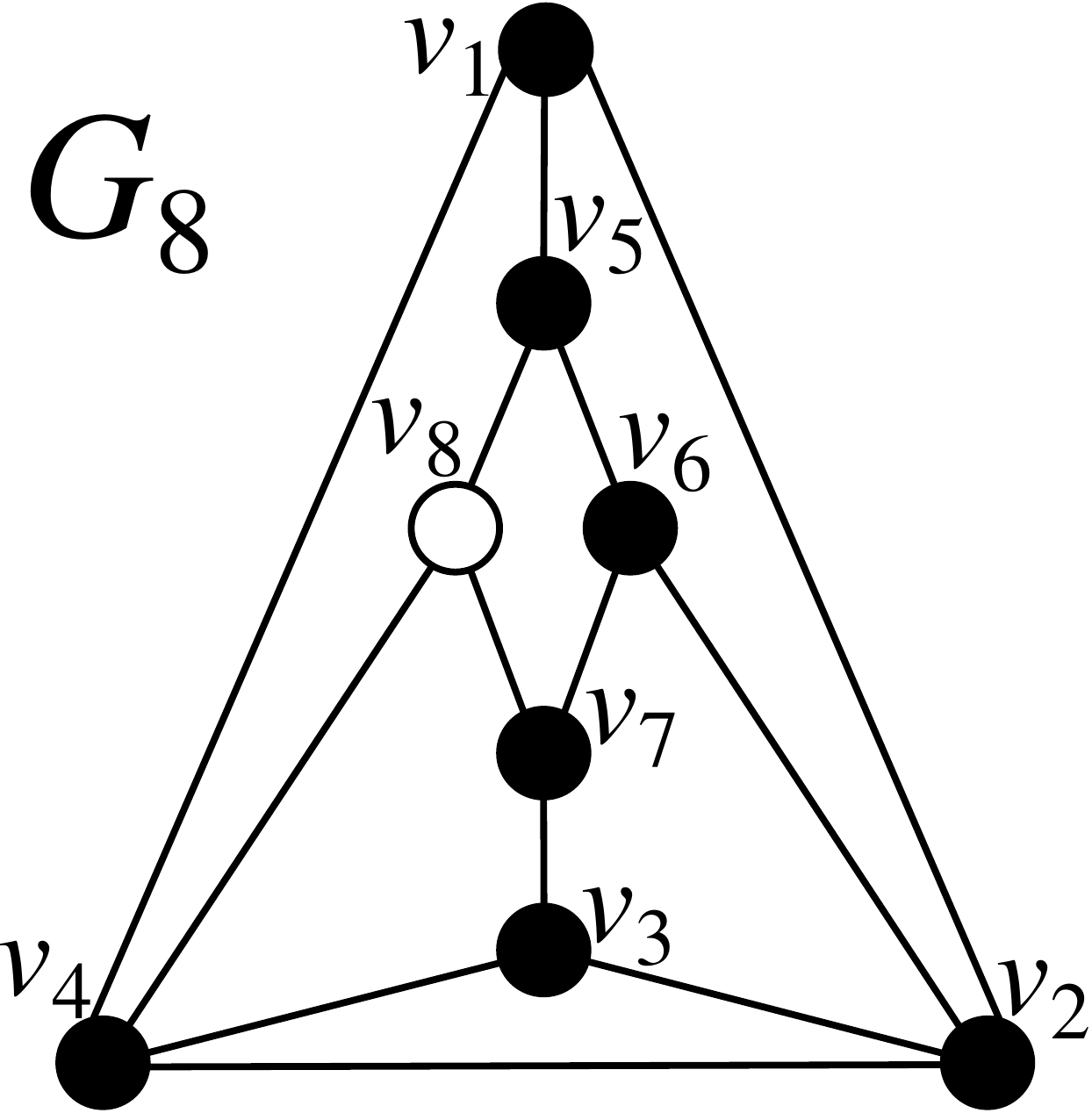} & \includegraphics[width=0.19\textwidth]{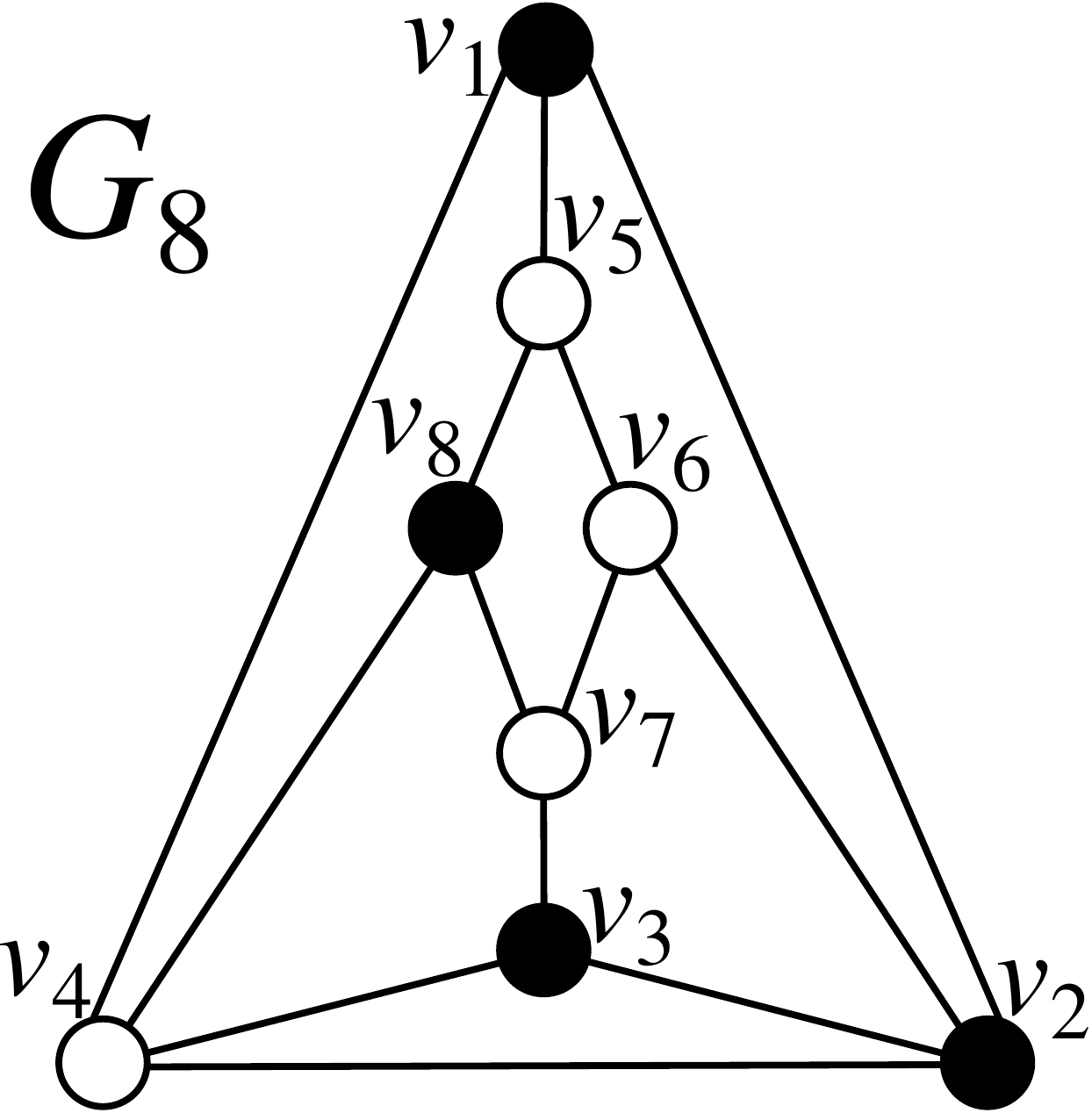} & \includegraphics[width=0.19\textwidth]{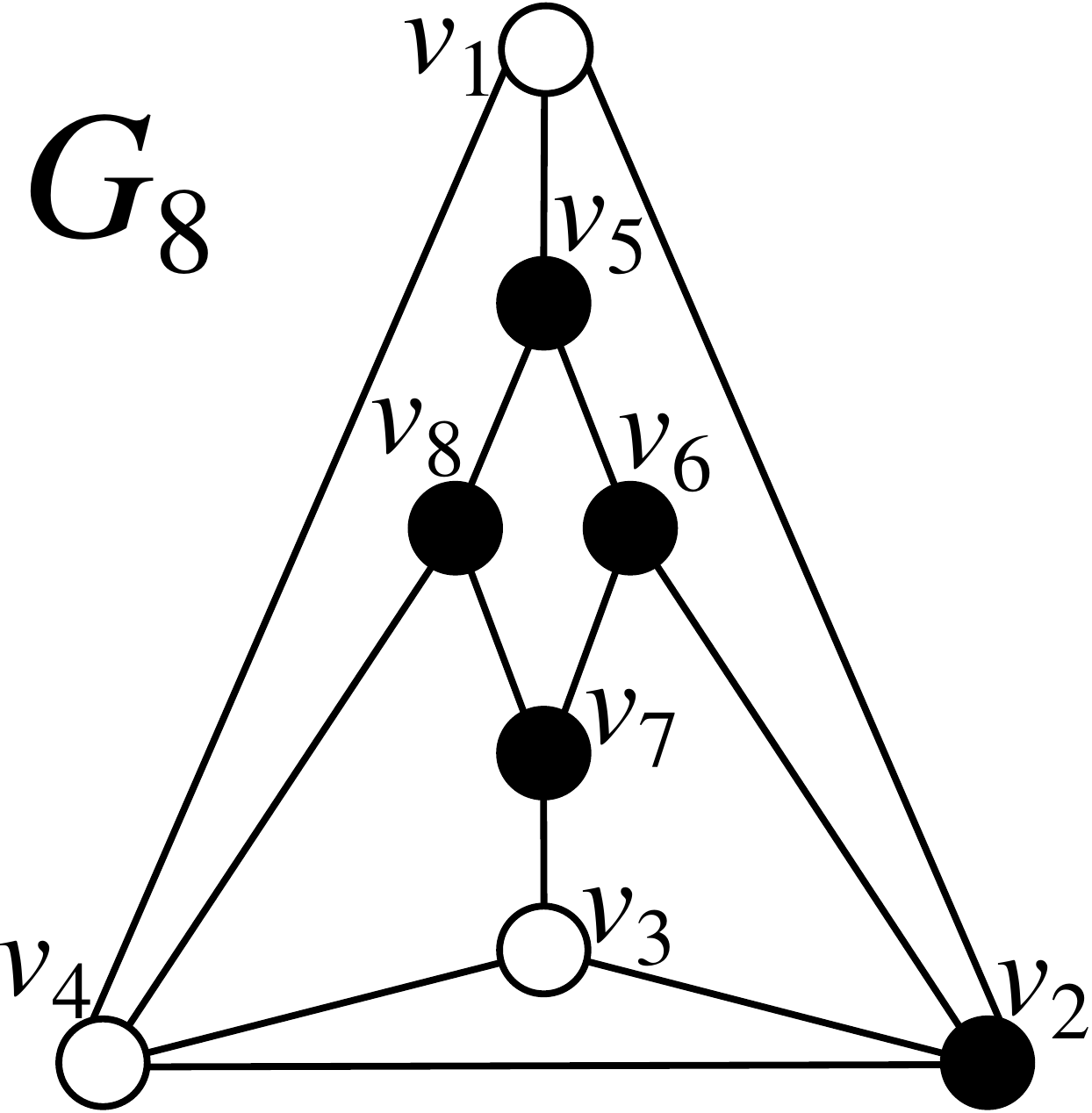} \\ (a) & (b) & (c) & (d)
    \end{tabular}
    \caption{Optimal OLD (a), RED:OLD (b), IC (c) and RED:IC (d) sets on $G_8$. Shaded vertices represent detectors}
    \label{fig:ex-finite}
\end{figure}

Figure~\ref{fig:ex-finite} shows examples of an OLD set, a RED:OLD set, an IC, and a RED:IC on $G_8$.
One can verify that both (a) and (c) satisfy Theorem~\ref{theo:det-sys} with $R_{det}(v)$ being $N(v)$ and $N[v]$, respectively; so (a) is an OLD set and (c) is an IC.
Similarly, both (b) and (d) satisfy Theorem~\ref{theo:red-det-sys} with $k=1$ and $R_{det}(v)$ being $N(v)$ and $N[v]$, respectively; so (b) is a RED:OLD set and (d) is a RED:IC.
Additionally, one can verify that for any of these sets, no smaller set satisfying the same requirements exists; therefore, these are optimal.
Thus, we see that OLD($G_8$) $= 5$, RED:OLD($G_8$) $= 7$, IC($G_8$) $= 4$, and RED:IC($G_8$) $= 5$.
If we would prefer to use densities, we also have that OLD\%($G_8$) $= \frac{5}{8}$, RED:OLD\%($G_8$) $= \frac{7}{8}$, IC\%($G_8$) $= \frac{4}{8}=\frac{1}{2}$, and RED:IC\%($G_8$) $= \frac{5}{8}$.

\section{Detection Systems on Infinite Graphs}

\begin{figure}[ht]
    \centering
    \begin{tabular}{cc}
        \begin{tabular}{c}\includegraphics[width=0.4\textwidth]{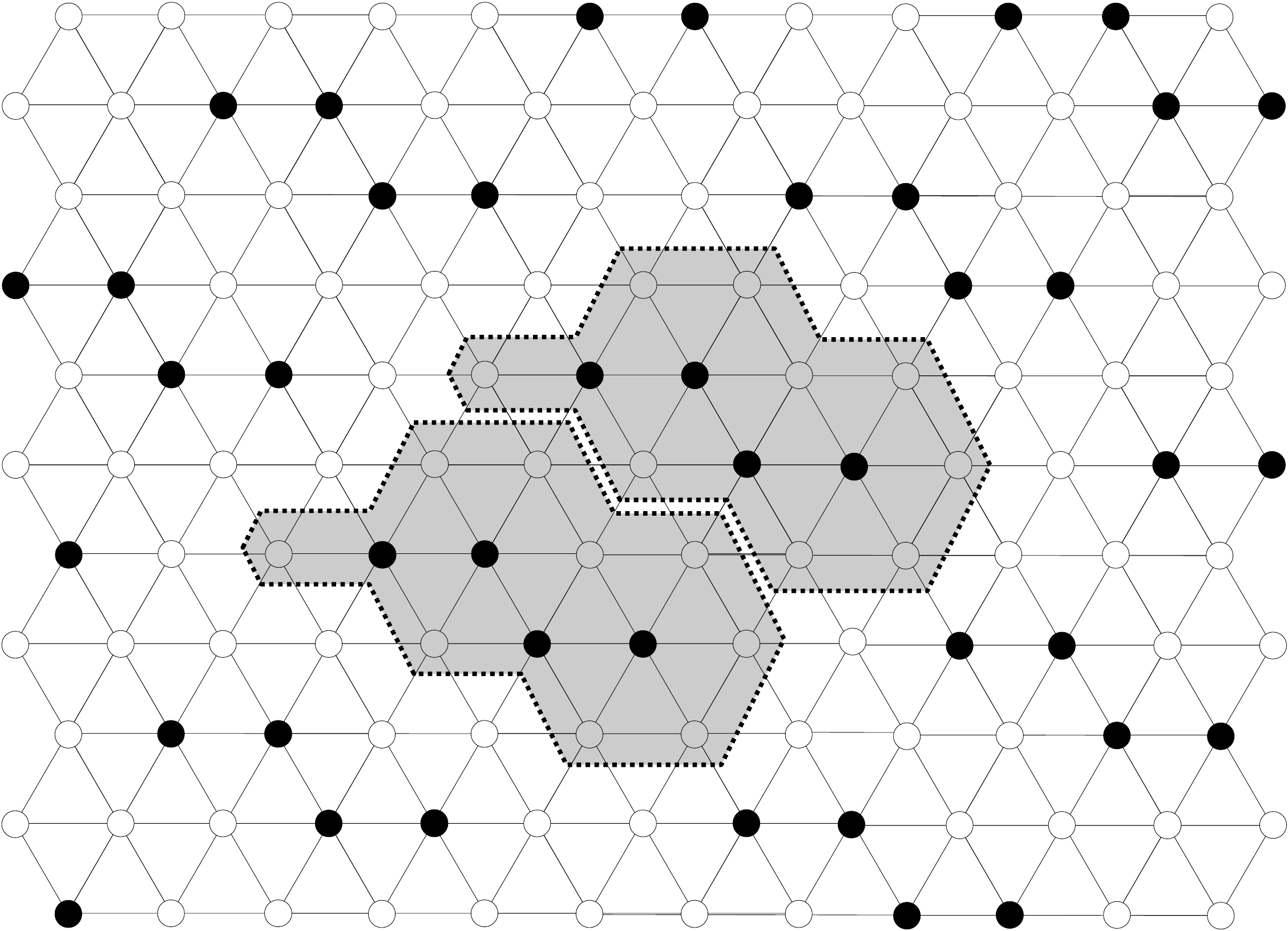}\\(a)\end{tabular} &
        \begin{tabular}{c}\includegraphics[width=0.4\textwidth]{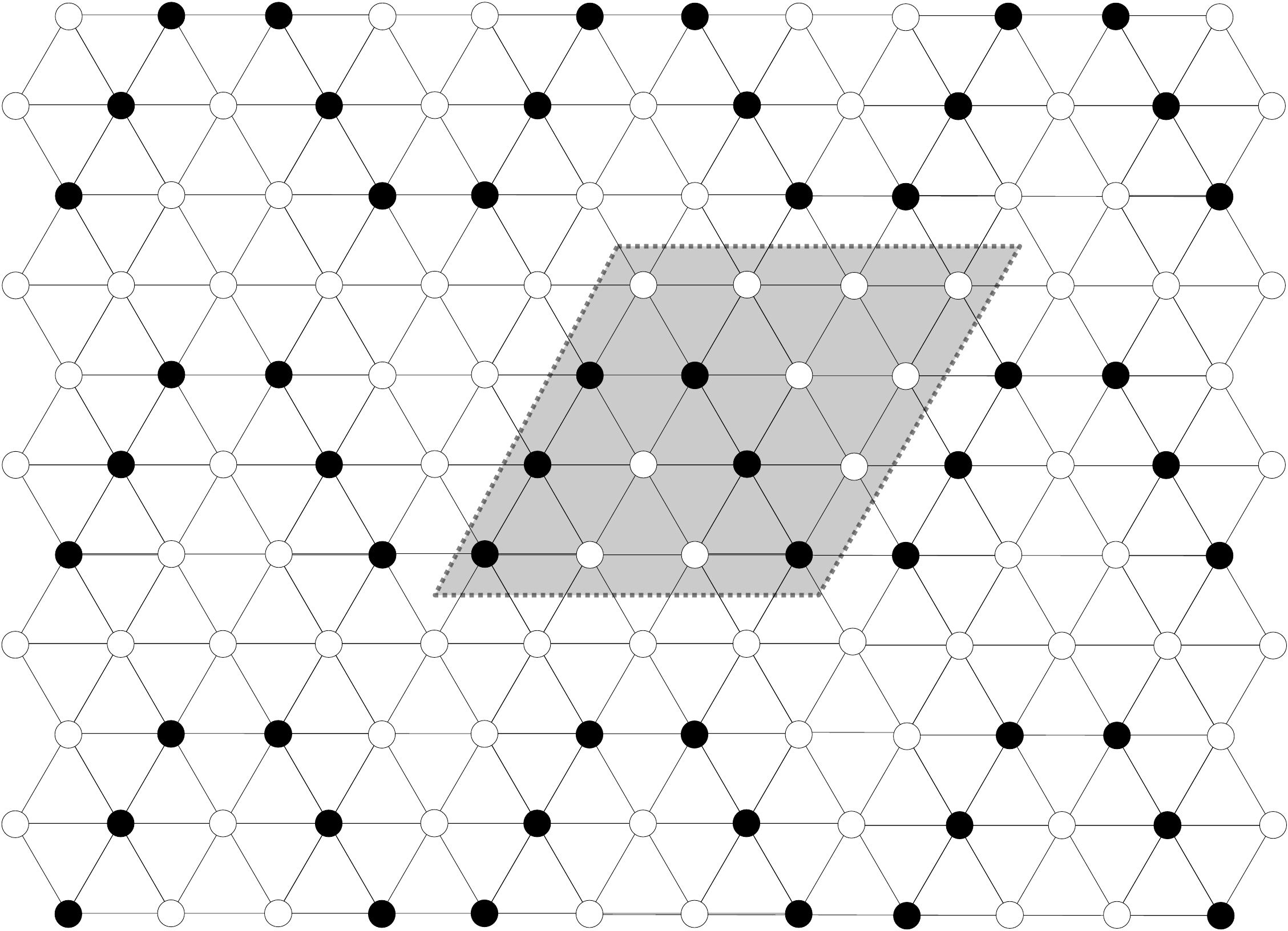}\\(b)\end{tabular}
    \end{tabular}
    \caption{Optimal OLD (a) and RED:OLD (b) sets on the TRI graph, constructed via tessellation. \cite{ourtri}}
    \label{fig:triexample}
\end{figure}

One common method for constructing vertex subsets for an infinite graph $G$ is via \emph{tessellation}, in which we take a non-empty region $R \subseteq V(G)$ of the graph and tile it to cover $V(G)$ with no overlapping of tiles.
Let $T$ be the set of all tiles required to cover $V(G)$.
$T$ must have the following properties: $R \in T$, $\cup_{A \in T}{A} = V(G)$, and $\forall A,B \in T$ with $A$$\neq$$B$, $|A|$$=$$|B| \;\land\; A \cap B$$=$$ \varnothing$.
Therefore, $T$ is a partition of $V(G)$.
As $|A| = |R|$ for any tile $A \in T$, there exists a bijective function $f_A: A \rightarrow R$ mapping from any vertex in tile $A$ back into the original tile $R$.
In tessellation, we consider every tile to be a copy of $R$, so $f_A$ maps any vertex in tile $A$ back to its corresponding geometric position in the original tile $R$.
Additionally, as $A \cap B = \varnothing$ for any two distinct tiles $A,B \in T$ and $\cup_{A \in T}{A} = V(G)$ there exists a surjective function $g: V(G) \rightarrow T$ mapping any vertex in the graph to its (only) enclosing local tile.
That is, $\forall v \in V(G),\; v \in g(v) \in T$.
Combining these functions, we can create a function $F: V(G) \rightarrow R$, $F(v) \equiv f_{g(v)}(v)$ mapping from any vertex in the graph back to its corresponding geometric position in the original tile $R$.
Once this tessellation has been constructed, we can use it to create the desired structured subsets of the infinite graph.
Let $S \subseteq R$ be the desired subset of the original tile $R$.
Then the tessellated subset generated by $S$, denoted $\Psi_S$, is $\cup_{v \in S}{F^{-1}(v)}$ where $F^{-1}$ is the total pre-image of $F$.
That is, $\forall v \in R,\; F^{-1}(v) \equiv \{w \in V(G) : F(w) = v\}$.

Figure~\ref{fig:triexample} shows two examples of tessellated constructions for an optimal OLD-set \cite{oldtri} and a RED:OLD set \cite{ftsets}, respectively, on the \emph{infinite triangular grid (TRI)}.
Let $R \subseteq V(G)$ be the original tile, and let $S \subseteq R$.
One of the benefits of tessellation is that, by their very construction, the density of any generated subset $\Psi_S$ in $V(G)$ is equal to the density of $S$ in $R$.
Thus, as Figure~\ref{fig:triexample} shows optimal solutions we see that OLD\%(TRI) $= \frac{4}{13}$ and RED:OLD\%(TRI) $= \frac{6}{16} = \frac{3}{8}$.

In this paper, we consider RED:OLD sets and RED:ICs on the \emph{infinite king's grid (K)}---also known as the \emph{infinite king's graph}---an 8-regular graph inspired by chess, with vertices representing positions on an infinite chessboard and edges representing all the normal movements of a king.
Several papers have been published exploring various graphical parameters on the king's grid.
The first such paper---by Cohen, Honkala, and Lobstein \cite{chlz}---covered the creation of ICs; importantly, they proved that $\frac{2}{9} \le$ IC\%($K$) $\frac{4}{17}$.
This was soon improved by Charon, Hudry, and Lobstein \cite{chl}, who demonstrated that IC\%($K$) $= \frac{2}{9}$.
Honkala and Laihonen \cite{hl2} were the first to study locating-dominating sets on $K$, proving that the minimum density is $\frac{1}{5}$.
Many papers have since been published covering various properties of ICs \cite{chhl,dantas,flp,hl3,p6} and LD sets \cite{p2,p4,p5} on finite and infinite king's grids.
We will explore upper and lower bounds for RED:OLD\%($K$) and RED:IC\%($K$).

\FloatBarrier
\section{Lower Bounds}\label{sec:lower-bound-common}

\begin{wrapfigure}{r}{0.25\textwidth}
    \centering
    \includegraphics[width=0.2\textwidth]{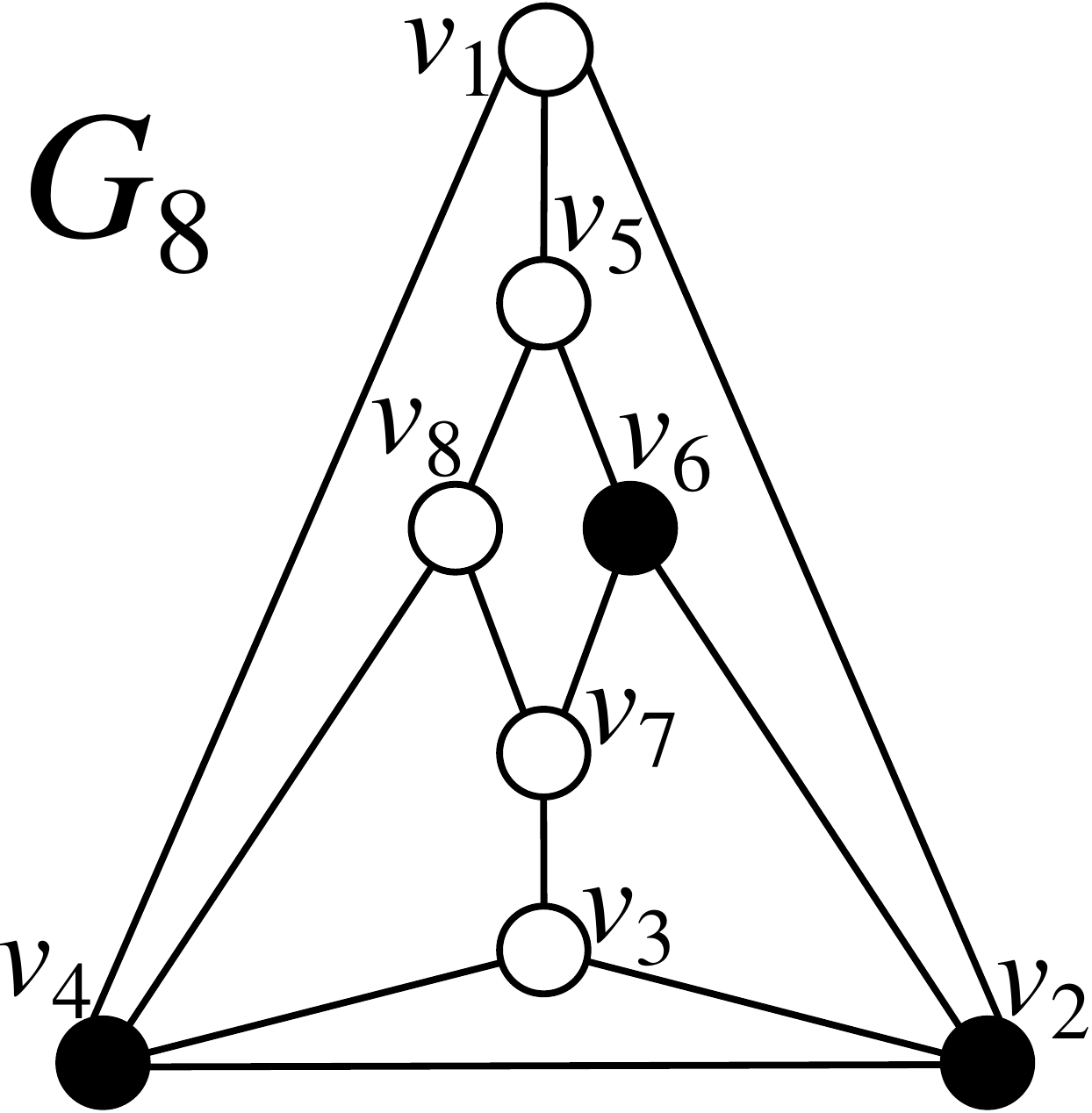}
    \caption{$G_8$ open-dominating set}
    \label{fig:ex-share}
\end{wrapfigure}

Depending on the graph, directly proving a lower bound for the minimum density of a detector-based distinguishing set can be difficult.
Instead, it's often easier to use what's known as a \emph{share argument}, established by Slater \cite{ftld} and used extensively by numerous papers on related concepts of distinguishing sets \cite{ourtri, ftsets, tmb, oldking, old}.
The inverse of average share is equal to the density.
Thus, if we find an upper bound for the average share, its inverse is a lower bound for density.

Let $G$ be a graph, $S \subseteq V(G)$ be a set of detectors, $v \in S$, and $u \in R_{det}(v)$.
Vertex $u$ is dominated exactly $dom_S(u)$ times, one of which is due to $v$.
Thus, $v$'s share in dominating $u$, denoted $sh_S[u](v)$, is $\frac{1}{dom_S(u)}$.
And the (total) share of $v$, denoted $sh_S(v)$, is $\sum_{w \in R_{det}(v)}{sh_S[w](v)}$.
We will also introduce the following compact notation: where $A \subseteq R_{det}(v)$, the \emph{partial share} of $v$ over $A$, denoted $sh_S[A](v)$, is $\sum_{w \in A}{sh_S[w](v)}$.

As an example of calculating shares and partial shares, consider the open-dominating set $S = \{v_2,v_4,v_6\}$ on the $G_8$ graph in Figure~\ref{fig:ex-share}.
Consider $v_2 \in S$, which is adjacent to $\{v_1,v_3,v_4,v_6\}$.
$dom_S(v_4) = dom_S(v_6) = 1$, so $sh_S[v_4](v_2) = sh_S[v_6](v_2 = \frac{1}{1} = 1$.
$dom_S(v_1) = dom_S(v_3) = 2$, so $sh_S[v_1](v_2) = sh_S[v_3](v_2) = \frac{1}{2}$.
Therefore, $sh_S(v_2) = sh_S[v_1,v_3,v_4,v_6](v_2) = 1 + 1 + \frac{1}{2} + \frac{1}{2} = 3$.
Now consider $v_4 \in S$, which is adjacent to $\{v_1,v_2,v_3,v_8\}$.
$dom_S(v_1) = dom_S(v_2) = dom_S(v_3) = 2$ and $dom_S(v_8) = 1$, so $sh_S(v_4) = sh_S[v_1,v_2,v_3,v_8](v_4) = \frac{1}{2} + \frac{1}{2} + \frac{1}{2} + 1 = \frac{5}{2}$.
And finally, we have $v_6 \in S$ with $sh_S(v_6) = \frac{1}{2} + 1 + 1 = \frac{5}{2}$.
The average share is $\frac{1}{3}\left[3 + \frac{5}{2} + \frac{5}{2}\right] = \frac{8}{3}$ and indeed we see that the density of this solution is $\frac{3}{8}$.

\begin{wrapfigure}{r}{0.25\textwidth}
    \centering
    \includegraphics[width=0.15\textwidth]{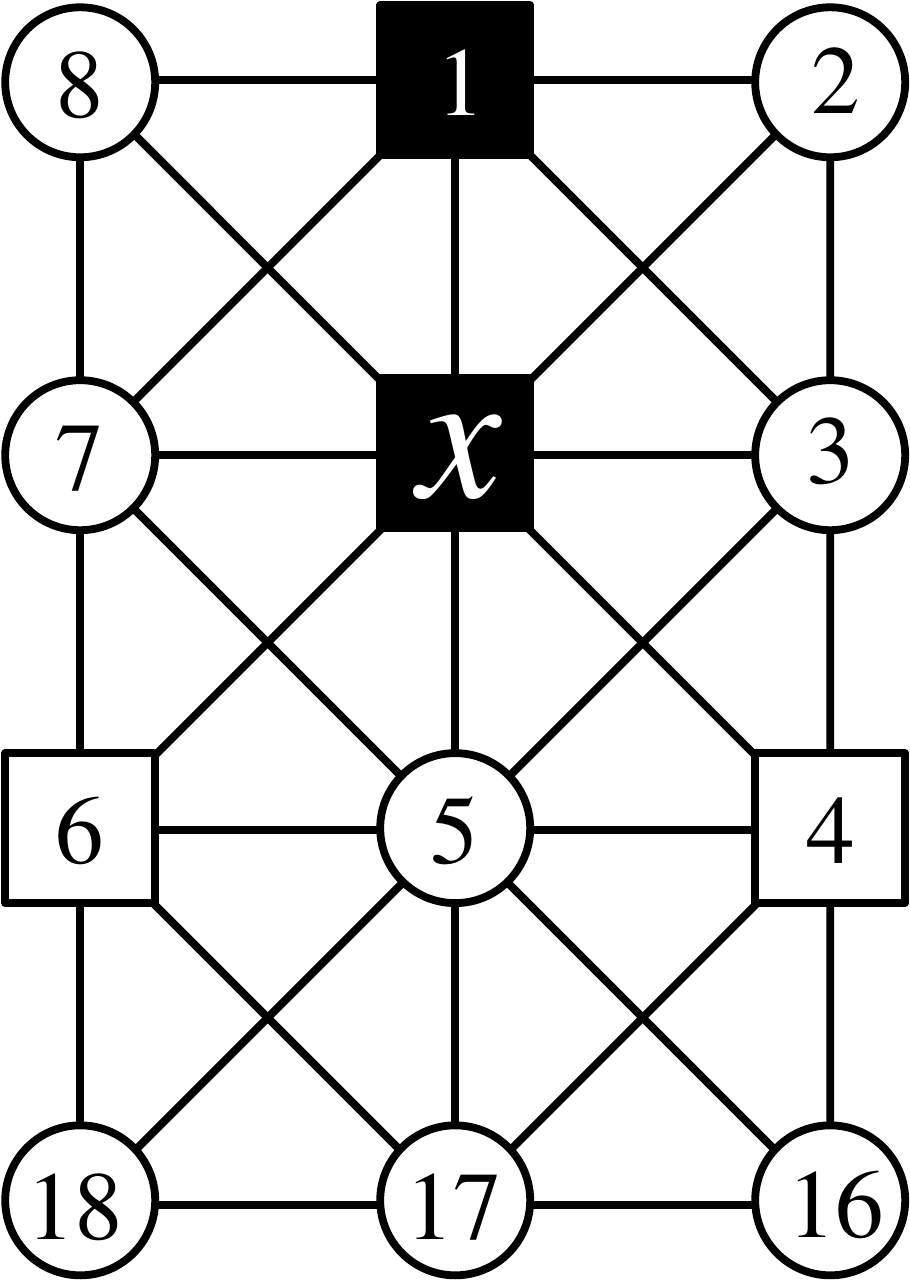}
    \caption{Example of drawing conventions}
    \label{fig:labeling-example}
\end{wrapfigure}
The graphs shown in subsequent figures use circles to represent ``unknowns" (vertices which may or may not be detectors), shaded squares to represent detectors, and non-shaded squares to represent non-detectors.
Additionally, any vertex labeled simply with number $k$ will be referred to as $v_k$ in the text.
For instance, in Figure~\ref{fig:labeling-example}, $\{x,v_1\}$ are detectors, $\{v_4,v_6\}$ are non-detectors, and everything else is currently unknown.

Let $S \subseteq V(K)$ be a RED:OLD set on $K$.
Consider the scenario in Figure~\ref{fig:labeling-example}.
Vertices $v_2$, $v_3$, $v_7$, and $v_8$ are dominated by the same two vertices.
If one of the four vertices is 2-dominated then the other three must be at least 4-dominated in order to be distinguished; thus, $sh_S[v_2,v_3,v_7,v_8](x) \le \frac{1}{2} + \frac{3}{4} = \frac{5}{4}$.
Otherwise all four vertices are at least 3-dominated, in which case $sh_S[v_2,v_3,v_7,v_8](x) \le \frac{4}{3}$.
In either case we see that $sh_S[v_2,v_3,v_7,v_8](x) \le \frac{4}{3}$.
This finding is useful, but it can easily be generalized to any graph and many possible configurations and types of sets:

\begin{lemma}\label{lem:common}
Suppose $G$ is a graph, $S \subseteq V(G)$ is a $k$-distinguishing set of detectors with uniform and symmetric detection regions $R_{det}$, $k \ge 0$, $D \subseteq S$ is a non-empty set of detectors with $v \in D$, and $A \subseteq V(G)$ is a non-empty set of vertices with $D \subseteq \cap_{u \in A}{R_{det}(u)}$.
If $\exists u \in A$ such that $dom_S(u) = |D|$ then $sh_S[A](v) \leq \frac{1}{|D|} + \frac{|A|-1}{|D|+k}$.
Otherwise, $sh_S[A](v) \leq \frac{|A|}{|D|+1}$.
In either case,
\begin{equation*}
sh_S[A](v) \le \begin{cases}
    \frac{|A|}{|D|+1}, & |A||D|(k-1) \ge (|D|+1)k \\
    \frac{1}{|D|} + \frac{|A|-1}{|D|+k}, & |A||D|(k-1) < (|D|+1)k
    \end{cases}
\end{equation*}
\end{lemma}
\begin{proof}
By hypothesis, $D \subseteq S$ and $D \subseteq \cap_{u \in A}{R_{det}(u)}$.
Because $R_{det}$ is uniform and symmetric, by Theorem~\ref{theo:loc-code}, $D \subseteq \cap_{u \in A}{\mathscr{L}_S(u)}$.
Suppose $\exists u \in A$ such that $dom_S(u) = |D|$.
Then $\mathscr{L}_S(u) = D$ and we observe that $\forall w \in A-\{u\},\; |\mathscr{L}_S(w) \triangle \mathscr{L}_S(u)| = |\mathscr{L}_S(w) - D| = |\mathscr{L}_S(w)| - |D| = dom_S(w) - |D|$.
By hypothesis, $S$ is $k$-distinguishing, so this quantity must be at least $k$.
Therefore, $dom_S(w)-|D| \ge k$, which of course means $dom_S(w) \ge |D|+k$.
Hence, $sh_S[A](v) \le \frac{1}{|D|} + \frac{|A|-1}{|D|+k}$.
Otherwise $\forall u \in A,\; dom_S(u) \ge |D|+1$, in which case $sh_S[A](v) \le \frac{|A|}{|D|+1}$.

\newpage
In either case $sh_S[A](v) \le \max\left\{\frac{|A|}{|D|+1}, \frac{1}{|D|} + \frac{|A|-1}{|D|+k}\right\}$.
We will now characterize when this is $\frac{|A|}{|D|+1}$.
For convenience, let $a = |A|$ and $d = |D|$:
\begin{align*}
    \frac{a}{d+1} &\ge \frac{1}{d} + \frac{a-1}{d+k} \\
    ad(d+k) &\ge (d+1)(d+k) + (a-1)(d+1)d \\
    ad^2 + adk &\ge d^2 + dk + d + k + ad^2 + ad - d^2 - d \\
    adk &\ge dk + k + ad\\
    ad(k-1) &\ge (d+1)k
\end{align*}
\end{proof}

\FloatBarrier
\section{Bounds for $\textrm{RED:OLD\%}(K)$}\label{sec:red-old}

\FloatBarrier
\subsection{Upper Bound}
\begin{figure}[ht]
    \centering
    \begin{tabular}{cc}
        \includegraphics[width=0.44\textwidth]{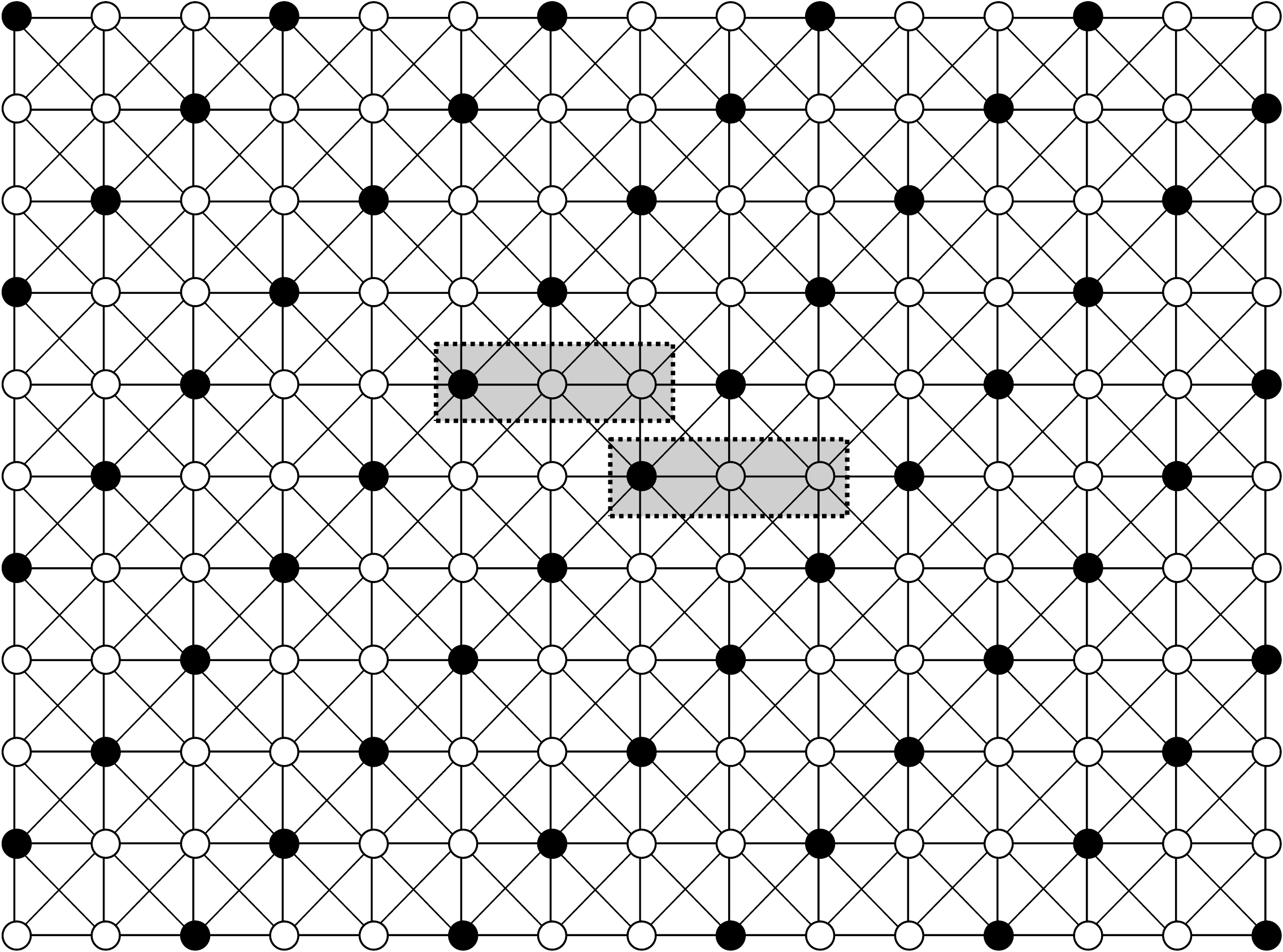} & \includegraphics[width=0.44\textwidth]{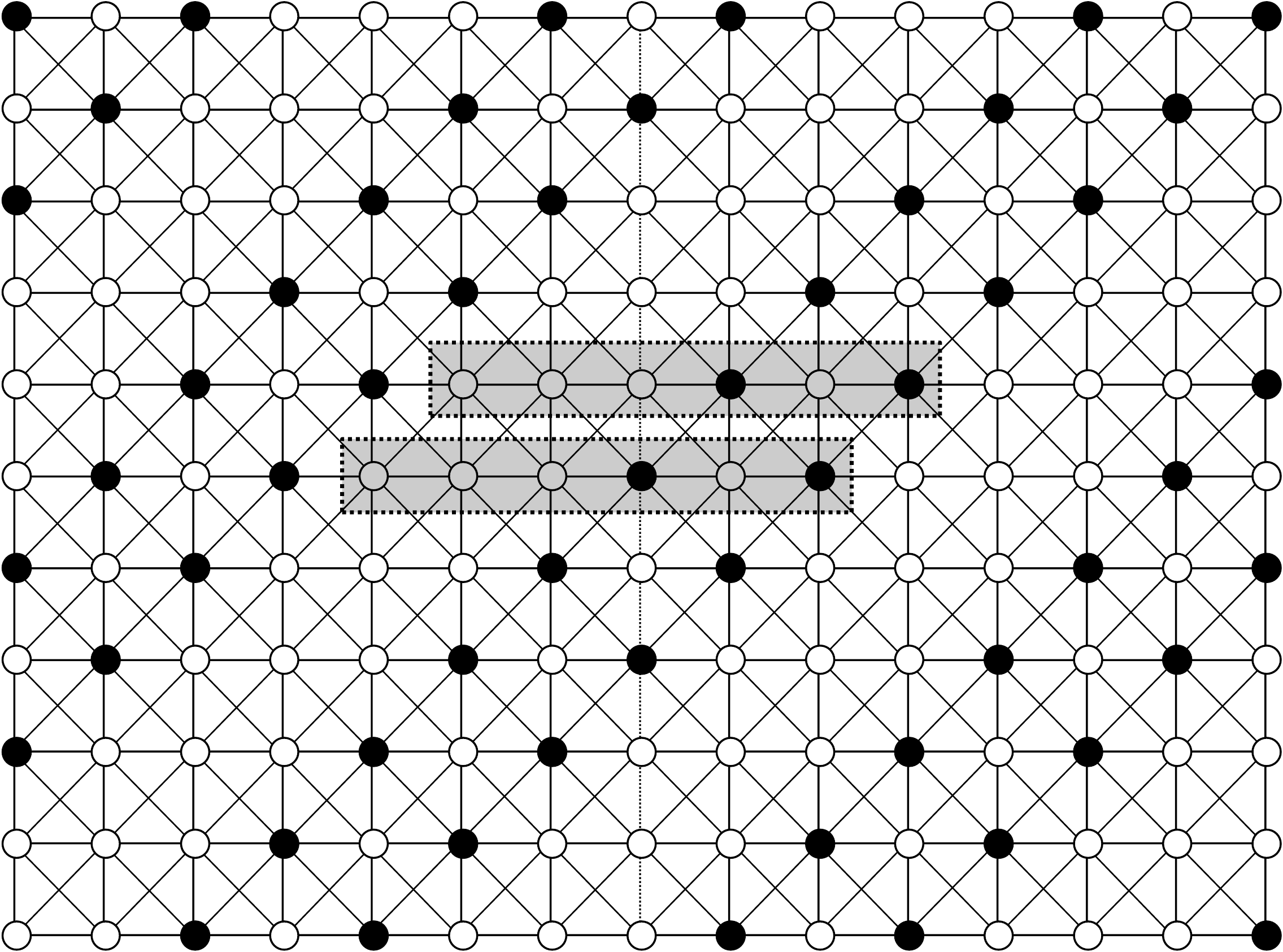} \\ (a) & (b)
    \end{tabular}
    \caption{Two non-isomorphic RED:OLD sets showing a density of $\frac{1}{3}$. Shaded vertices are in the RED:OLD set}
    \label{fig:soln}
\end{figure}

We have found two non-isomorphic RED:OLD sets on the infinite king's grid $K$, as shown in Figure~\ref{fig:soln}.
These solutions can be formed via tessellation of rectangles (regions outlined in the figure).
For Figure~\ref{fig:soln}~(a), we see that 1 of the 3 vertices in the tile is a detector (shaded vertices), giving a density of $\frac{1}{3}$.
Similarly, in Figure~\ref{fig:soln}~(b) we see that 2 of the 6 vertices in the tile are used, also giving a density of $\frac{2}{6} = \frac{1}{3}$.
It is also possible to form these solutions via tessellation of diagonal ribbons of infinite length.
For the first solution (a), the ribbon consists of 1 diagonal of detectors and 2 diagonals of non-detectors, giving a density of $\frac{1}{3}$.
For the second solution (b), the ribbon consists of 2 diagonals of detectors and f diagonals of non-detectors, giving a density of $\frac{2}{6} = \frac{1}{3}$.

As we have a RED:OLD set on $K$ with density $\frac{1}{3}$, we have an upper bound for the optimal density: RED:OLD\%($K$) $\le \frac{1}{3}$.

\FloatBarrier
\subsection{Lower Bound}
Let $S \subseteq V(K)$ be a RED:OLD set on $K$, and let the notations $sh$, $dom$, and $k$-dominated implicitly use $S$.
These will persist for the rest of Section~\thesubsection.

$S$ is a 1-redundant detection system, so Theorem~\ref{theo:red-det-sys} gives us that every vertex must be at least 2-dominated, so for any detector $x \in S$ and any neighbor $v \in N(x)$, $sh[v](x) \le \frac{1}{2}$.
Additionally, $K$ is an 8-regular graph, so $|N(x)| = 8$.
Therefore, $sh(x) \le \frac{8}{2} = 4$, giving us a simple lower bound of RED:OLD\%($K$) $\ge \frac{1}{4}$.

\begin{wrapfigure}{r}{0.3\textwidth}
    \centering
    \includegraphics[width=0.25\textwidth]{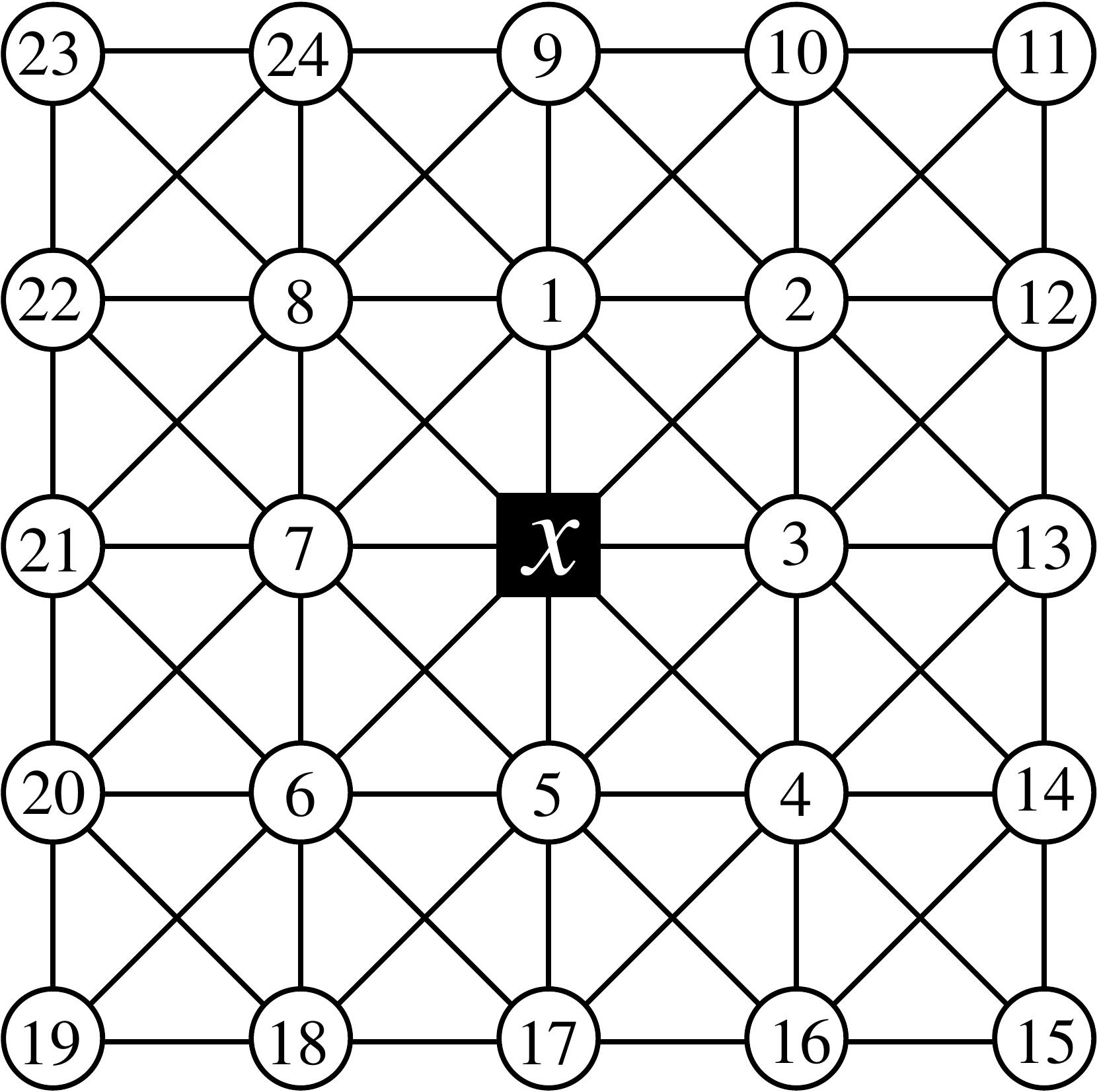}
    \caption{Simple bound example with labels}
    \label{fig:simple-bound}
\end{wrapfigure}
We will now prove another, higher bound.
Consider a detector $x \in S$ and its as-of-yet unknown local region of $K$, as shown in Figure~\ref{fig:simple-bound}.
Suppose that $\forall w \in N(x),\; dom(w) = 2$.
All neighbors of $x$ are dominated by $x$ itself.
Because the neighbors of $x$ are all 2-dominated by hypothesis, they must all be dominated exactly one more time.
However, if two neighbors of $x$ are dominated by the same two detectors, then they will not be distinguished.
If $v_9 \in S \lor v_{10} \in S$ then $\{v_1,v_2\}$ are not distinguished.
Likewise, if $v_1 \in S$ then $\{v_2,v_3\}$ are not distinguished.
And if $v_2 \in S$ then $\{v_1,v_3\}$ are not distinguished.
By symmetry, the only vertices in the figure which can possibly be detectors are $\{v_{11},v_{15},v_{19},v_{23}\}$, but this means that $\forall w \in \{v_1,v_3,v_5,v_7\},\; dom(w) = 1$, contradicting that $S$ is a RED:OLD set on $K$.
Therefore, $\exists w \in N(x)$ such that $dom(w) \ge 3$, so $sh[w](x) \le \frac{1}{3}$.
By Theorem~\ref{theo:red-det-sys}, any RED:OLD set at least 2-dominates every vertex, thus, where $u \in N(x)$, $dom(u) \ge 2$, thus $sh[u](x) \le \frac{1}{2}$.
Thus, $sh_S(x) \le \frac{7}{2} + \frac{1}{3} = \frac{23}{6}$, meaning we have a higher lower bound: RED:OLD\%($K$) $\ge \frac{6}{23}$.

At this point, we have two lower bounds for RED:OLD\%($K$): $\frac{1}{4}$ and $\frac{6}{23}$.
These are both still relatively far from our upper bound of $\frac{1}{3}$.
To help bridge the gap between upper and lower bounds, we would now like to prove a significantly more ambitious lower bound of $\frac{3}{10}$, which will require demonstrating an upper bound for average share of $\frac{10}{3}$.

\begin{theorem}\label{theo:main-sh}
In any RED:OLD set $S$ on $K$, the average share over all vertices in $S$ is no larger than $\frac{10}{3}$.
\end{theorem}

\cbeginproof
Let $x \in S$ be an arbitrary vertex in the RED:OLD set $S$.
We will explore all non-isomorphic configurations and show that the average share it gives rise to is at most $\frac{10}{3}$.

\newpage
\FloatBarrier
\newcase{1}
\begin{wrapfigure}{r}{0.2\textwidth}
    \vspace{-2.5em}
    \centering
    \includegraphics[width=0.15\textwidth]{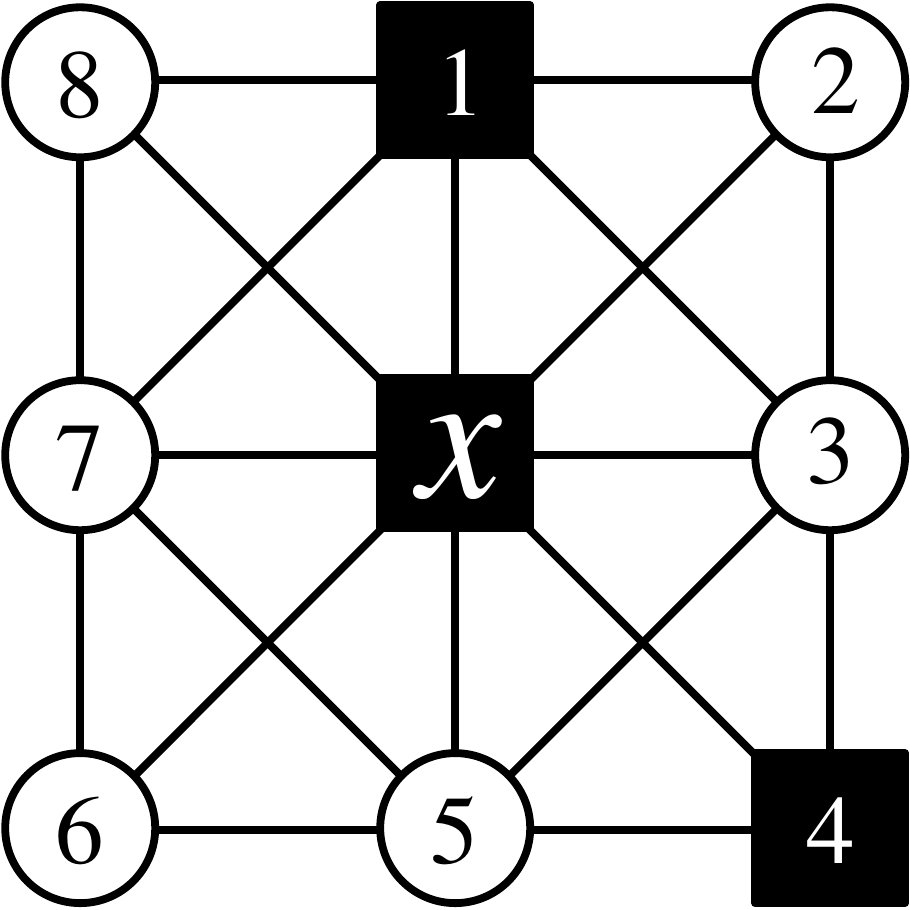}
    \caption{Case 1}
    \label{fig:red-old-case-1}
\end{wrapfigure}

We begin by considering when $x$ is dominated by two neighbors in a ``bent" shape, as in Figure~\ref{fig:red-old-case-1}.
We will demonstrate that in any possible sub-configuration $sh(x) \le \frac{13}{4} < \frac{10}{3}$.

From the figure we clearly see that $dom(v_3) \ge 3$.
Suppose $dom(v_3) = 3$.
By Lemma~\ref{lem:common}, $sh[v_2,v_3,v_7,v_8](x) \le \frac{4}{3}$.
Additionally, $dom(v_5) \ge 3$ in order to be distinguished from $v_3$.
Since each of the other three vertices in $N(x)$ is at least 2-dominated, we have $sh[v_1,v_4,v_6](x) \le \frac{3}{2}$.
Thus, $sh(x) \le \frac{4}{3} + \frac{1}{3} + \frac{3}{2} = \frac{19}{6} < \frac{13}{4}$.
Otherwise $dom(v_3) \ge 4$.
Then, by Lemma~\ref{lem:common}, $sh[v_2,v_7,v_8](x) \le 1$, and each of the other four vertices in $N(x)$ is at least 2-dominated, we have $sh[v_1,v_4,v_5,v_6](x) \le \frac{4}{2}$.
Thus, $sh(x) \le 1 + \frac{1}{4} + \frac{4}{2} = \frac{13}{4}$.

\FloatBarrier
\newcase{2}
\begin{wrapfigure}{r}{0.2\textwidth}
    \vspace{-2.5em}
    \centering
    \includegraphics[width=0.15\textwidth]{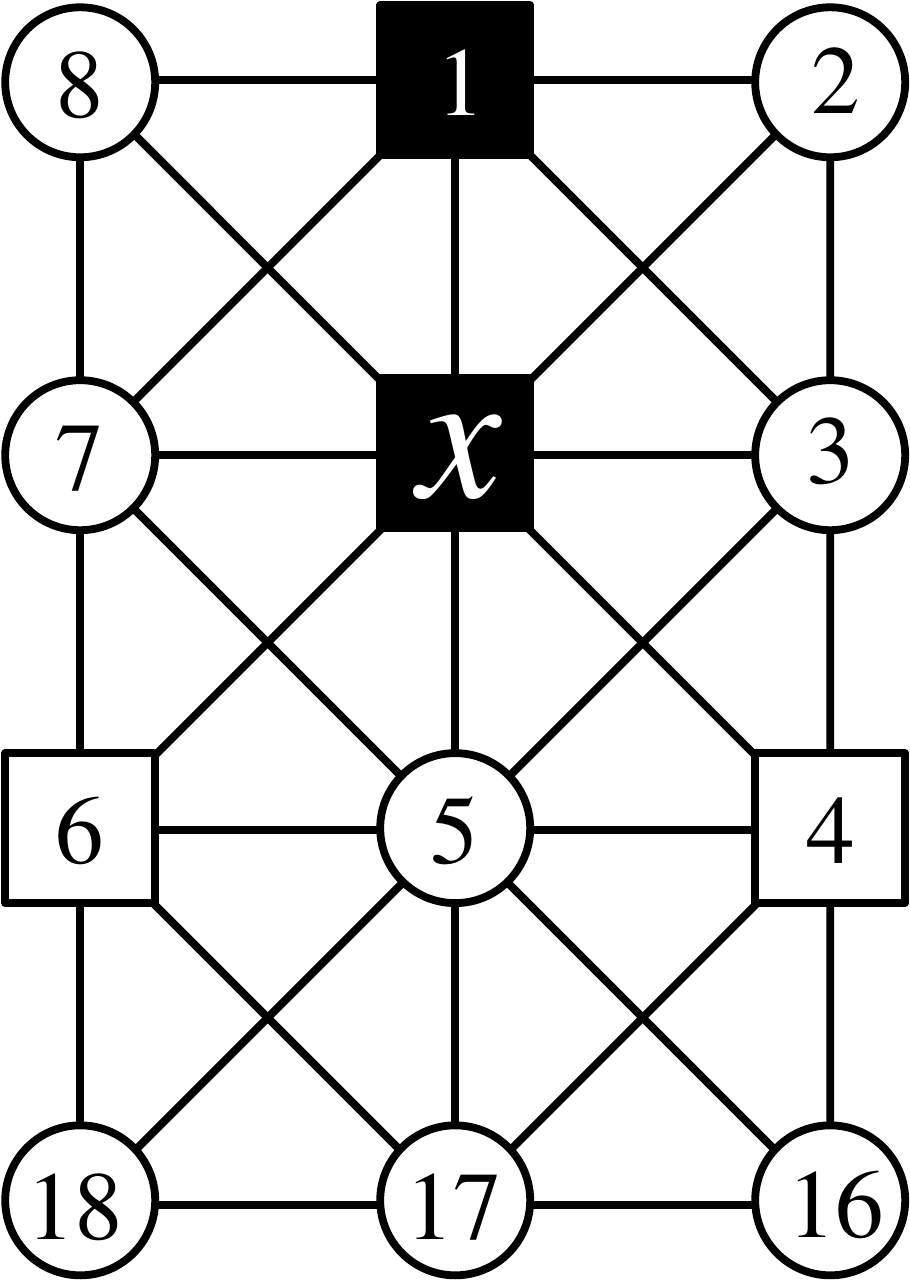}
    \caption{Case 2}
    \label{fig:red-old-case-2}
\end{wrapfigure}

Next, we consider when $x$ is dominated by vertex $v_1$ directly above it, as in Figure~\ref{fig:red-old-case-2}.
We will demonstrate that in any possible sub-configuration $sh(x) \le \frac{13}{4} < \frac{10}{3}$.

From Case~1 we know that if $\{v_4,v_6\} \cap S \neq \varnothing$ then $sh(x) \le \frac{13}{4}$ and we are done.
Therefore, we assume that neither of these is a detector.
Suppose that $\forall w \in \{v_4,v_5,v_6\},\; dom(w) = 2$.
If $v_{17} \in S$ then $\{v_4,v_5\}$ are impossible to distinguish due to both being 2-dominated by hypothesis; thus, $v_{17} \notin S$.
By the same logic, $v_{16} \notin S$ and $v_3 \notin S$.
Likewise, if $v_{18} \in S$ then $\{v_5,v_6\}$ are impossible to distinguish.
By the same logic, $v_7 \notin S$.
We therefore find that $dom(v_5) = 1$, contradicting that $S$ is a RED:OLD set.
Therefore, $\exists w \in \{v_4,v_5,v_6\}$ such that $dom(w) \ge 3$.
Additionally, by Lemma~\ref{lem:common}, $sh[v_2,v_3,v_7,v_8](x) \le \frac{4}{3}$.
Therefore, $sh(x) \le \frac{4}{3} + \frac{1}{3} + \frac{3}{2} = \frac{19}{6} < \frac{13}{4}$.

\FloatBarrier
\begin{wrapfigure}[5]{r}{0.2\textwidth}
    \vspace{-1em}
    \centering
    \includegraphics[width=0.15\textwidth]{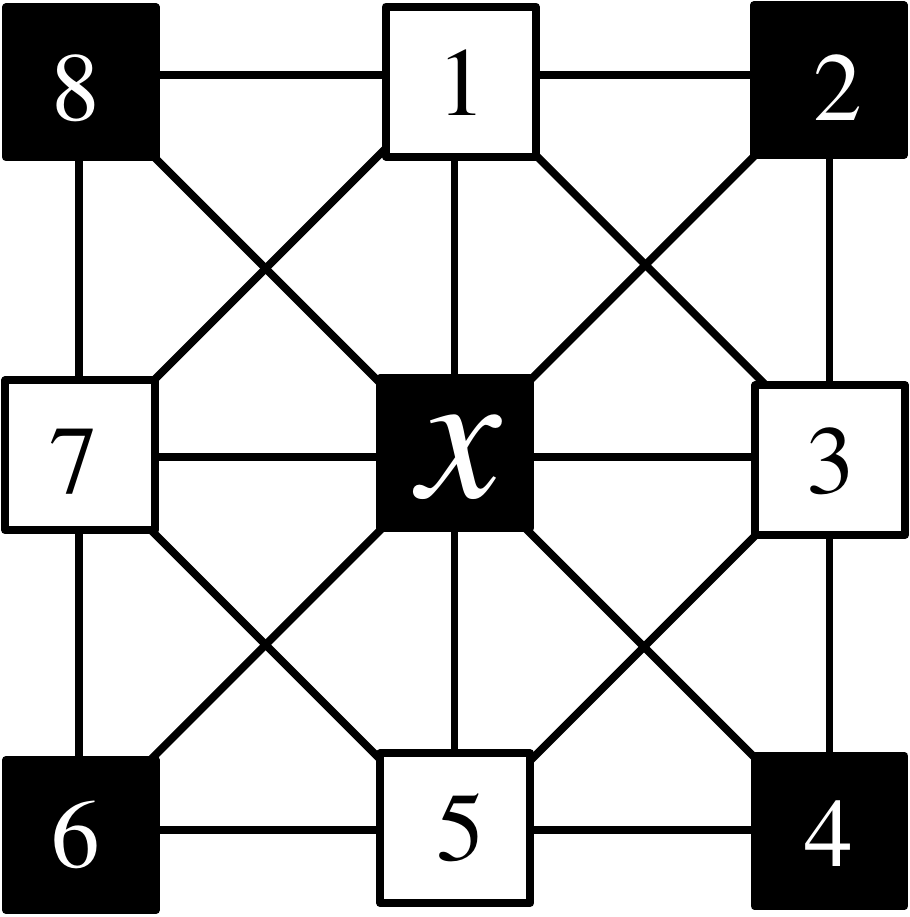}
    \caption{Case 3}
    \label{fig:red-old-case-3}
\end{wrapfigure}
\newcase{3}

We now consider when $x$ is dominated by four neighbors in an ``x" shape, as shown in Figure~\ref{fig:red-old-case-3}.
From Cases~1--2 we see that if $\{v_1,v_3,v_5,v_7\} \cap S \neq \varnothing$ then $sh(x) \le \frac{13}{4}$ and we are done.
Otherwise, we observe that $\forall w \in \{v_1,v_3,v_5,v_7\},\; dom(w) \ge 3$.
Thus $sh(x) \le \frac{4}{3} + \frac{4}{2} = \frac{10}{3}$.

\FloatBarrier
\newcase{4}

\begin{wrapfigure}{r}{0.3\textwidth}
    \centering
    \includegraphics[width=0.25\textwidth]{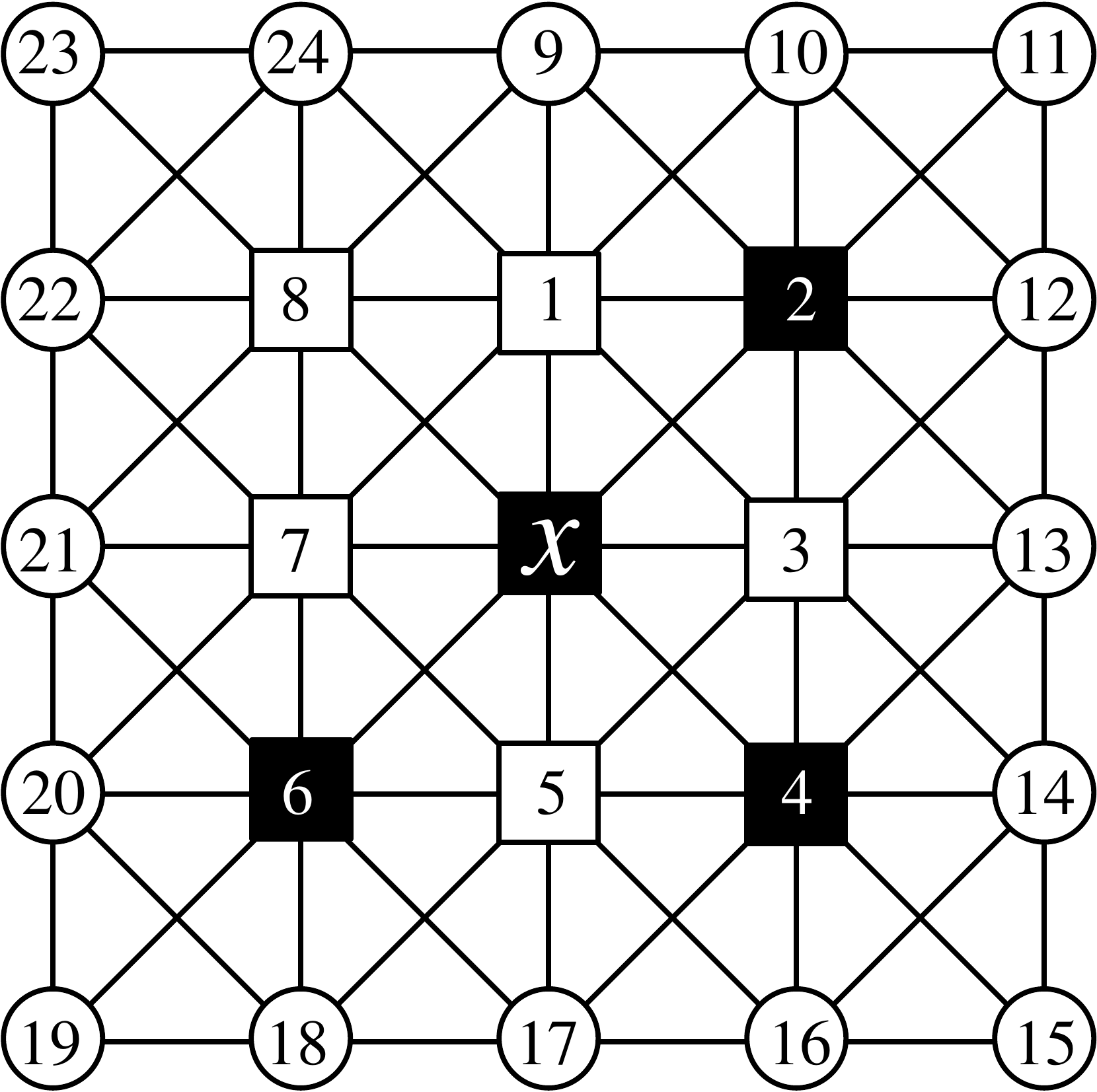}
    \caption{Case 4}
    \label{fig:red-old-case-4}
\end{wrapfigure}

We now consider the case where $x$ is dominated by three neighbors in a ``lambda" shape, as in Figure~\ref{fig:red-old-case-4}.
From Cases~1--3 we know that if $\{v_1,v_3,v_5,v_7,v_8\} \cap S \neq \varnothing$ then $sh(x) \le \frac{10}{3}$.
Therefore we assume $\{v_1,v_3,v_5,v_7,v_8\} \cap S = \varnothing$ and consider the sub-cases depending the values of $dom(v_3)$ and $dom(v_5)$ as follows.

First, suppose $dom(v_3) = 6$ or $dom(v_5) = 6$. If $dom(v_3) = 6$, then $\{v_{12},v_{13},v_{14}\} \subseteq S$.
This makes $dom(v_2) \ge 3$ and $dom(v_4) \ge 3$.
Additionally, we have that $dom(v_5) \ge 3$ by hypothesis.
Therefore, $sh(x) \le \frac{1}{6} + \frac{3}{3} + \frac{4}{2} = \frac{19}{6} < \frac{10}{3}$.
By symmetry, we get the same results when $dom(v_5) = 6$.

Next, we consider when $dom(v_3) = 3$.  We observe that $dom(v_3) = 3$ requires  $dom(v_1) \ge 3$.
If we have $dom(v_5) = 3$, then it also requires $dom(v_7) \ge 3$, hence we obtain $sh(x) \le \frac{4}{3} + \frac{4}{2} = \frac{10}{3}$ and we are done.
Therefore, we assume that $dom(v_5) \ge 4$.
Currently, we have three vertices, namely $v_1$, $v_3$, and $v_5$, which are at least 3-dominated. If $\exists w \in \{v_2,v_4,v_6,v_7,v_8\}$ such that $dom(w) \ge 3$ then $sh(x) \le \frac{4}{3} + \frac{4}{2} = \frac{10}{3}$ and are done.
Therefore, we assume $\forall w \in \{v_2,v_4,v_6,v_7,v_8\},\; dom(w) = 2$. 
By hypothesis, we have $dom(v_7) = 2$, thus $\{v_{20},v_{21},v_{22}\} \cap S = \varnothing$.
Presently, $dom(v_8) = 1$ and there are three unknown vertices, $\{v_9,v_{23},v_{24}\}$ which can possibly second-dominate $v_8$.
If $v_9 \in S$, then $\{v_2,v_8\}$ are not distinguished, a contradiction.
If $v_{23} \in S$ then $v_{24} \notin S $ leaving only $v_{10}$ to have $dom(v_1) \ge 3$, but this makes $\{v_1,v_2\}$ impossible to distinguish, a contradiction.
Therefore, the only way to have $dom(v_8) = 2$ requires $v_{24} \in S$, and this results in $dom(v_1) \ge 4$ because $\{v_1,v_8\}$ would not be distinguished, otherwise.
Thus $sh(x) \le \frac{2}{4} + \frac{1}{3} + \frac{5}{2} = \frac{10}{3}$.
By symmetry, if $dom(v_5) = 3$, we get the same results.

Now, we will proceed with the remaining possibilities for the domination of $\{v_3,v_5\}$. Due to symmetry, without loss of generality let $4 \le dom(v_3) \le dom(v_5) \le 5$ and we consider the following three configurations; $dom(v_3) = 5 = dom(v_5)$, $dom(v_3) = 4 \le 5 = dom(v_5)$, and $dom(v_3) = 4 = dom(v_5)$.

If we want to achieve $dom(v_3) = 5 = dom(v_5)$, any possible way to do this results in $dom(v_4) \ge 3$.
Thus $sh(x) \le \frac{2}{5} + \frac{1}{3} + \frac{5}{2} = \frac{97}{30} < \frac{10}{3}$, $dom(v_3) = 4 \le 5 = dom(v_5)$, 

\begin{wrapfigure}[14]{r}{0.3\textwidth}
    \centering
    \includegraphics[width=0.25\textwidth]{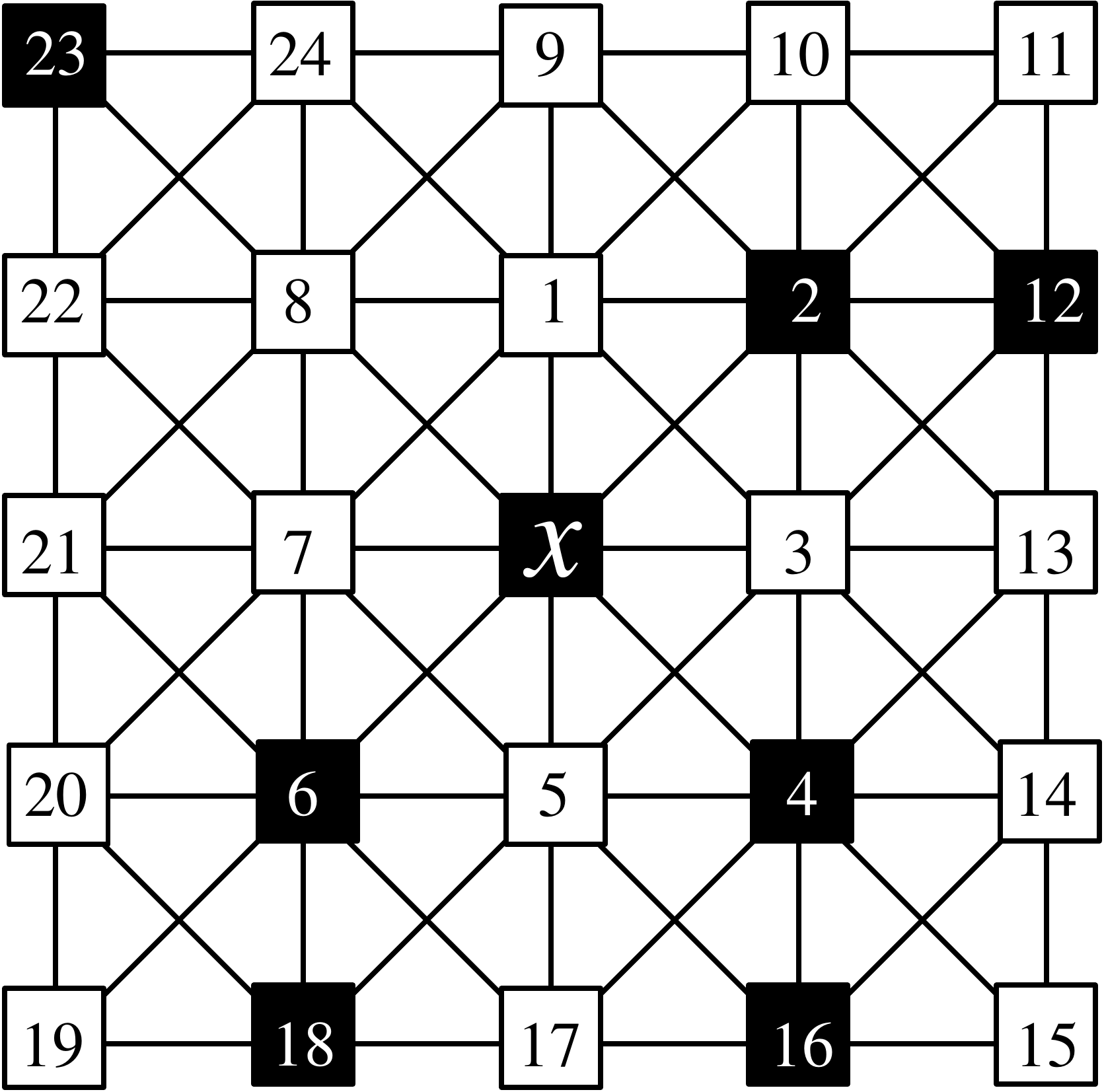}
    \caption{First terminal configuration of Case~4 with $sh(x) > \frac{10}{3}$}
    \label{fig:red-old-case-4-big-1}
\end{wrapfigure}

Now, we assume that $dom(v_3) = 4 \le 5 = dom(v_5)$.
We initially observe that if $\exists w \in \{v_1,v_2,v_4,v_6,v_7,v_8\}$ such that $dom(w) \ge 3$ then $sh(x) \le \frac{1}{5} + \frac{1}{4} + \frac{1}{3} + \frac{5}{2} = \frac{197}{60} < \frac{10}{3}$.
Thus, we need only consider when $\forall w \in \{v_1,v_2,v_4,v_6,v_7,v_8\},\; dom(w) = 2$.
Suppose $v_{17} \in S$, then $dom(v_4) \ge 3 \lor dom(v_6) \ge 3$ in order to be distinguished, a contradiction.
By hypothesis, $dom(v_5) = 5$, so $v_{17} \notin S \land \{v_{16},v_{18}\} \subseteq S$.
Because $dom(v_1) = 2$, $\{v_9,v_{10},v_{24}\} \cap S = \varnothing$.
Similarly, $dom(v_7) = 2$, requiring $\{v_{20},v_{21},v_{22}\} \cap S = \varnothing$.
Then $v_{23} \in S$ to make $dom(v_8) \ge 2$, and $v_{19} \notin S$ because $dom(v_6) = 2$.
We also see that $dom(v_4) = 2$, which means $\{v_{13},v_{14},v_{15}\} \cap S = \varnothing$.
This means $v_{12} \in S$ because $dom(v_3) = 4$, and $v_{11} \notin S$ because $dom(v_2) = 2$.
This leads to the sub-configuration shown in Figure~\ref{fig:red-old-case-4-big-1}.
We see that $sh(x) = \frac{1}{5} + \frac{1}{4} + \frac{6}{2} = \frac{69}{20} > \frac{10}{3}$.
However, we observe that $x$ is adjacent to $\{v_2,v_4,v_6\} \subseteq S$, all of which have ``bent" shape.
As demonstrated in Case~3, $\forall w \in \{v_2,v_4,v_6\},\; sh(w) \le \frac{13}{4} < \frac{10}{3}$.
Additionally, the only other detectors these 3 vertices are adjacent to are $\{v_{12},v_{16},v_{18}\}$, all of which are adjacent via vertical or horizontal edges.
Thus, the findings of Cases~1--4 tell us that $\forall w \in \{v_{12},v_{16},v_{18}\},\; sh(w) \le \frac{13}{4} < \frac{10}{3}$.
Therefore, we are safe to use $\{x,v_2,v_4,v_6\}$ in an averaging argument; consider the average share: $\frac{1}{4} \left[ \frac{69}{20} + \frac{13}{4} + \frac{13}{4} + \frac{13}{4} \right] = \frac{33}{10} < \frac{10}{3}$.
Hence, we have that the average share is still no larger than $\frac{10}{3}$.
\vspace{1.5em}

\begin{figure}[ht]
    \centering
    \begin{tabular}{c@{\hskip 5em}c}
        \includegraphics[width=0.25\textwidth]{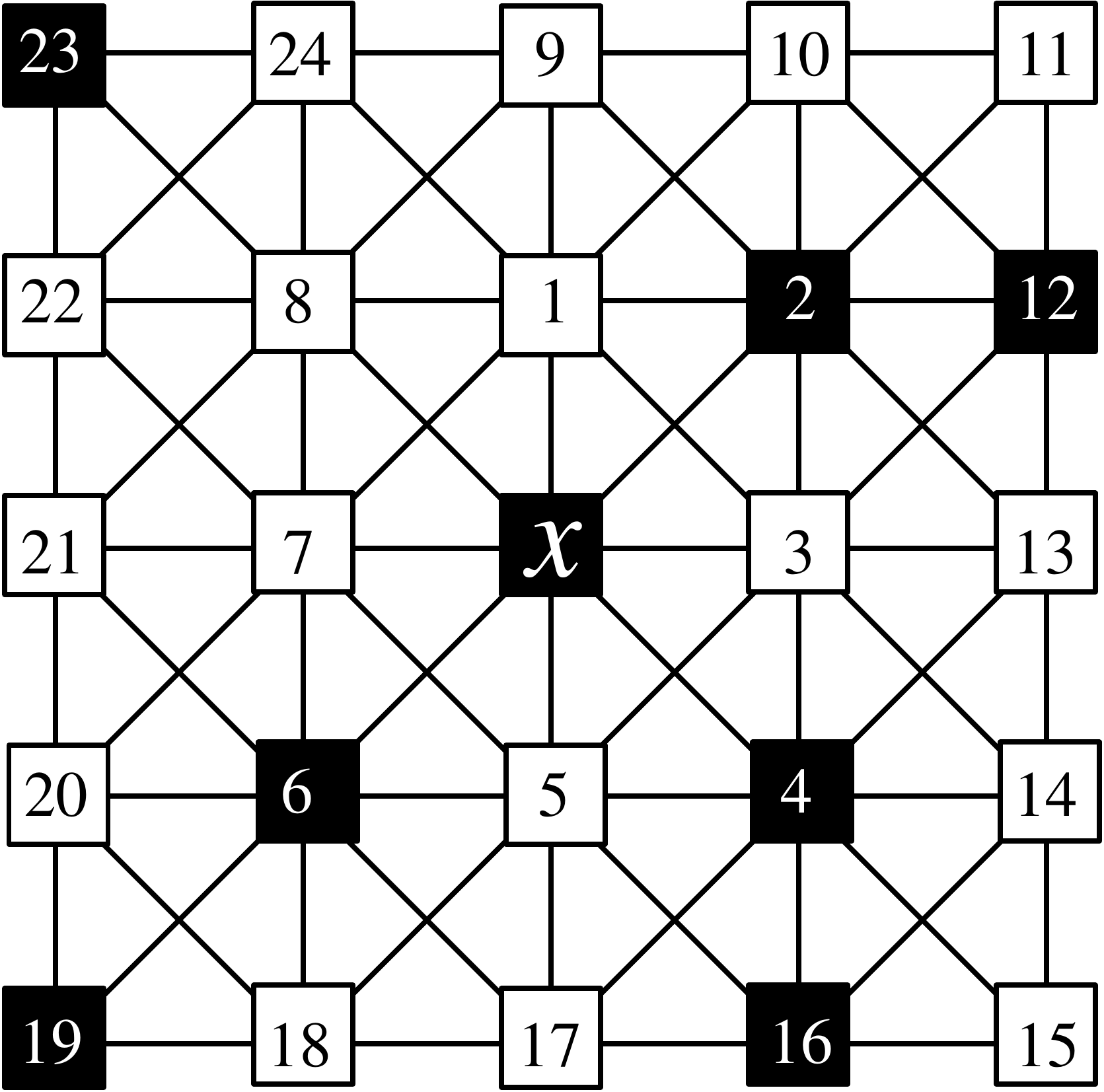} & \includegraphics[width=0.25\textwidth]{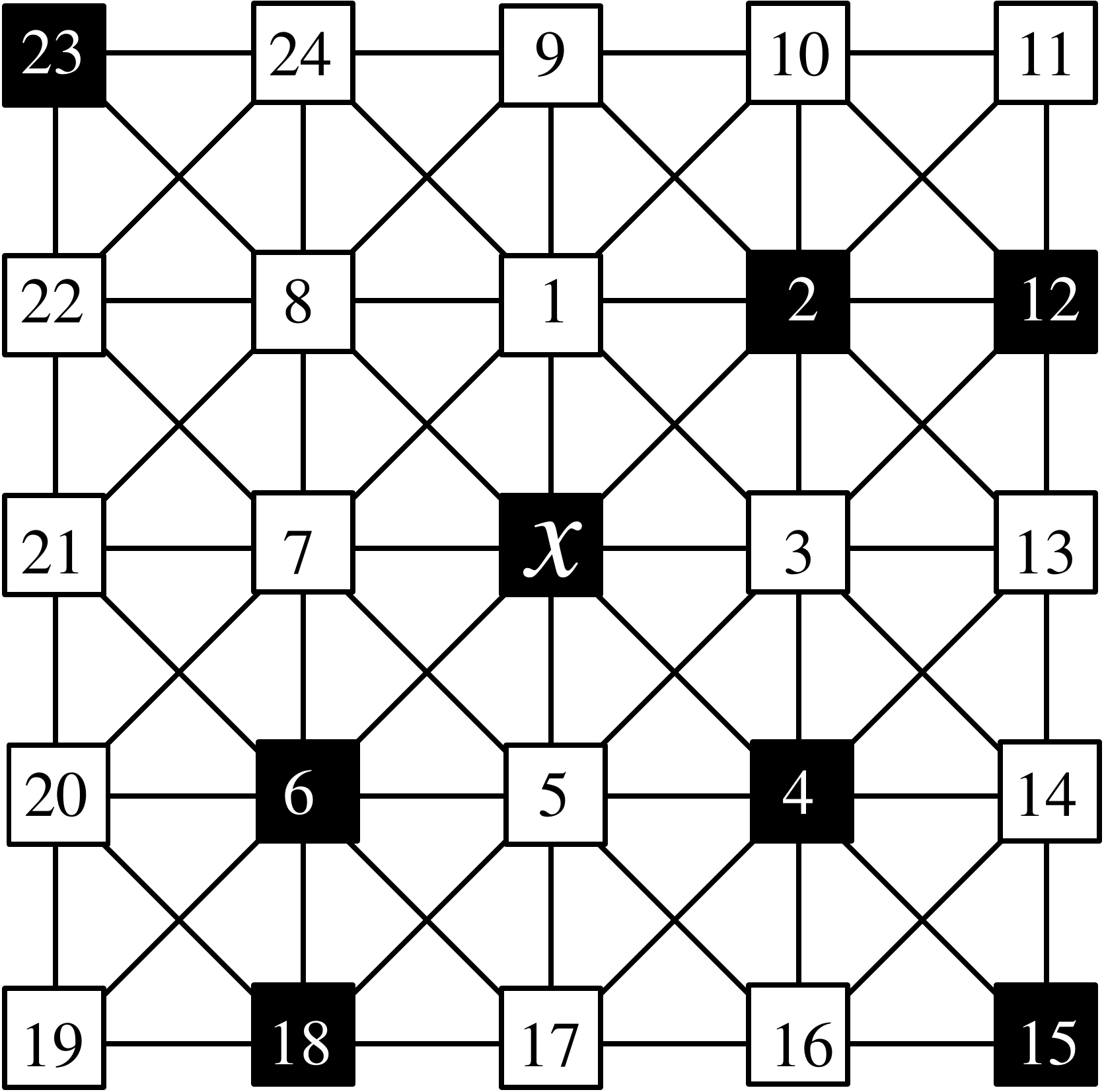} \\ (a) & (b)
    \end{tabular}
    \caption{Second and third terminal configurations of Case~4 with $sh(x) > \frac{10}{3}$}
    \label{fig:red-old-case-4-big-23}
\end{figure}

As the last sub-case, we have $dom(v_3) = 4 = dom(v_5)$.
We initially observe that if $\exists w \in \{v_1,v_2,v_4,v_6,v_7,v_8\}$ such that $dom(w) \ge 3$ then $sh(x) \le \frac{1}{5} + \frac{1}{4} + \frac{1}{3} + \frac{5}{2} = \frac{197}{60} < \frac{10}{3}$.
Thus, we need only consider when $\forall w \in \{v_1,v_2,v_4,v_6,v_7,v_8\},\; dom(w) = 2$.
Because $dom(v_4) = 2$, at most one of the additional detectors for $v_3$ and $v_5$ can be in $N(v_4)$.
Due to symmetry, without loss of generality let the additional detector for $v_3$ not be in $N(v_4)$; $\{v_{13},v_{14}\} \subseteq N(v_4)$, so $v_{12} \in S \land v_{13} \notin S \land v_{14} \notin S$.
Additionally, because $dom(v_2) = 2$ by hypothesis, $v_{11} \notin S$.
We also see that $dom(v_1) = dom(v_7) = 2$, implying $\{v_9,v_{10},v_{20},v_{21},v_{22},v_{24}\} \cap S = \varnothing$.
This requires $v_{23} \in S$ to make $dom(v_8) \ge 2$.
If $v_{17} \in S$ then $\{v_4,v_6\}$ are impossible to distinguish; thus, $v_{17} \notin S$.
By hypothesis, $dom(v_6) = 2$; thus exactly one of $\{v_{18},v_{19}\}$ must be a detector.
We will now consider two cases: $v_{19} \in S$ or $v_{19} \notin S$.

If $v_{19} \in S$ then $v_{18} \notin S$.
By hypothesis, $dom(v_5) = 5$, so $v_{16} \in S$, and by hypothesis $dom(v_4) = 2$, requiring $v_{15} \notin S$.
This leads to the sub-configuration shown in Figure~\ref{fig:red-old-case-4-big-23}~(a).
$sh(x) = \frac{2}{4} + \frac{6}{2} = \frac{7}{2} > \frac{10}{3}$.
However, as demonstrated previously in this case we can safely average with the adjacent ``bent" shape vertices $\{v_2,v_4\}$, which have share no more than $\frac{13}{4} < \frac{10}{3}$.
$\frac{1}{3} \left[ \frac{7}{2} + \frac{13}{4} + \frac{13}{4} \right] = \frac{10}{3}$.
Thus in this terminal sub-configuration, we still have that the average share is no more than $\frac{10}{3}$.

If $v_{19} \notin S$ then $v_{18} \in S$.
Then, by similar logic to the previous paragraph, $v_{16} \notin S \land v_{15} \in S$.
We now have the sub-configuration shown in Figure~\ref{fig:red-old-case-4-big-23}~(b).
This again has $sh(x) = \frac{2}{4} + \frac{6}{2} = \frac{7}{2} > \frac{10}{3}$.
However, we also still have 2 ``bent" shape neighbors of $x$: $\{v_2,v_6\}$, which are safe to use for averaging by the same logic, as demonstrated previously in this case.
Thus, we once again have the average share $\frac{1}{3} \left[ \frac{7}{2} + \frac{13}{4} + \frac{13}{4} \right] = \frac{10}{3}$.

\FloatBarrier
\newcase{5}
\begin{wrapfigure}{r}{0.3\textwidth}
    \centering
    \includegraphics[width=0.25\textwidth]{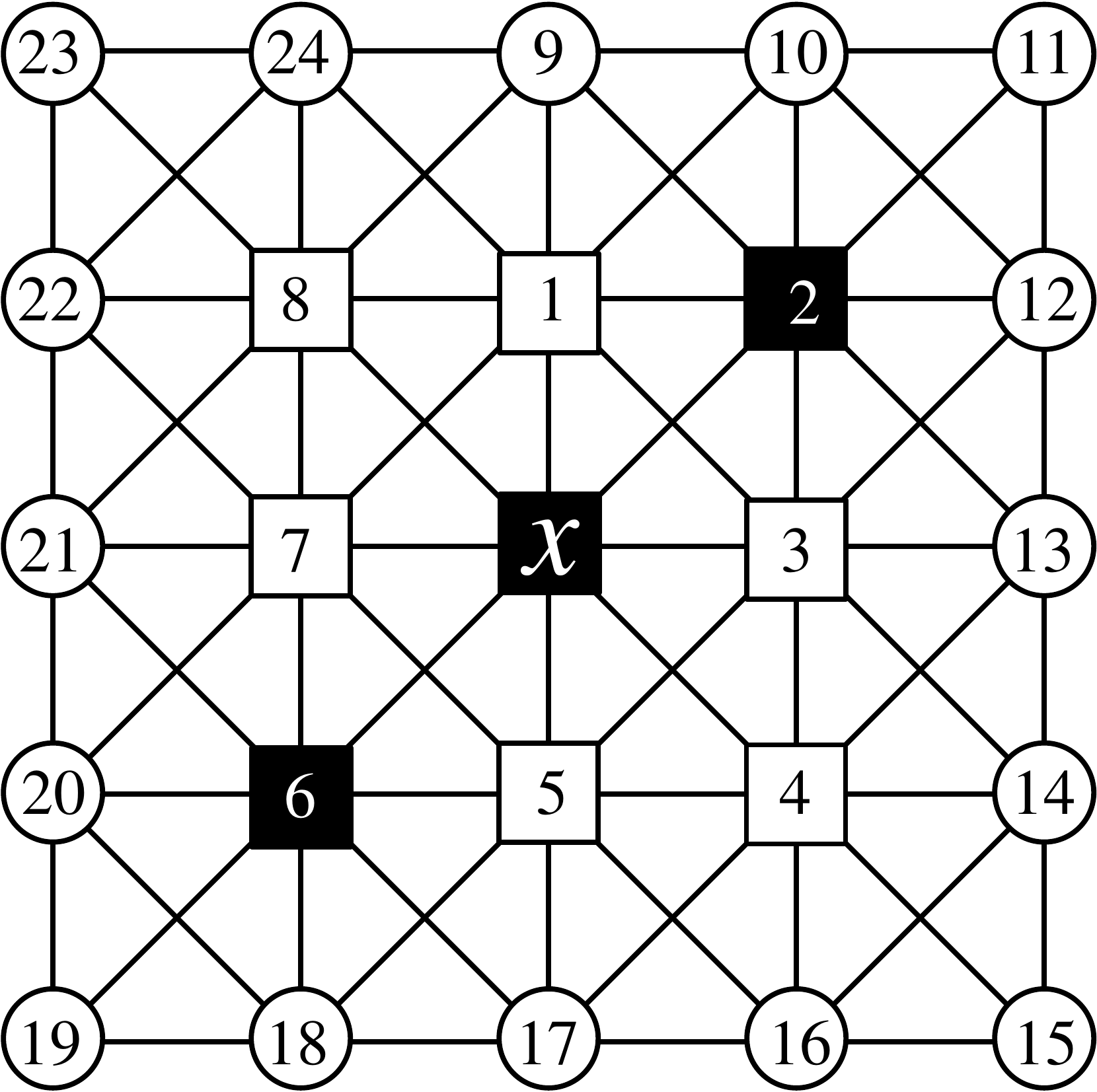}
    \caption{Case 5}
    \label{fig:red-old-case-5}
\end{wrapfigure}

We now consider when $x$ is dominated by two neighbors in a diagonal, as in Figure~\ref{fig:red-old-case-5}.
By applying the results of Cases~1--4 we see that if $\{v_1,v_3,v_4,v_5,v_7,v_8\} \cap S \neq \varnothing$ then the average share is no larger than $\frac{10}{3}$ and we are done.
Thus we need only consider when $\{v_1,v_3,v_4,v_5,v_7,v_8\} \cap S = \varnothing$.
We observe that $dom(v_1) \le 5$.
Suppose that $dom(v_1) = 5$, then $dom(v_2) \ge 3 \land dom(v_8) \ge 3$.
Additionally, $\{v_5,v_7\}$ are not distinguished, so at least one of them must be at least 3-dominated.
Thus $sh(x) \le \frac{1}{5} + \frac{3}{3} + \frac{4}{2} = \frac{16}{5} < \frac{10}{3}$.
By symmetry if $\exists w \in \{v_1,v_3,v_5,v_7\}$ such that $dom(w) \ge 5$ then $sh(x) \le \frac{10}{3}$.
Therefore, we need only consider when $\forall w \in \{v_1, v_3, v_5, v_7\},\; 2 \le dom(w) \le 4$.
We will break these into sub-cases by considering the vertex pair $(v_1, v_3)$.
Due to symmetry, without loss of generality, let $dom(v_1) \le dom(v_3)$.

If $dom(v_1) = 4 = dom(v_3)$, then $dom(v_2) \ge 3$.
Thus, $sh(x) \le \frac{2}{4} + \frac{1}{3} + \frac{5}{2} = \frac{10}{3}$.

Now consider when $dom(v_1) = 3 \le 4 = dom(v_3)$.
$\{v_5,v_7\}$ are not yet distinguished, so at least one of them is at least 3-dominated.
If the additional third detector for $v_1$ is $v_9$ or $v_{10}$ then $dom(v_2) = 3$; thus, $sh(x) \le \frac{3}{3} + \frac{1}{4} + \frac{4}{2} = \frac{13}{4} < \frac{10}{3}$.
So we need only consider when the additional third detector for $v_1$ is $v_{24} \in S$.
We then observe that $\{v_1,v_8\}$ are not distinguished, requiring $dom(v_8) \ge 3$.
Thus $sh(x) \le \frac{3}{3} + \frac{1}{4} + \frac{4}{2} = \frac{13}{4} < \frac{10}{3}$.

Next consider when $dom(v_1) = 3 = dom(v_3)$.
Again note that $\{v_5,v_7\}$ are not distinguished, so at least one of them must be at least 3-dominated.
We currently have three neighbors of $x$ which are at least 3-dominated; if we find one more, then $sh(x) \le \frac{4}{3} + \frac{4}{2} = \frac{10}{3}$.
If both of the additional detectors around $v_1$ and $v_3$ are in $N(v_2)$ then $dom(v_2) \ge 3$ and we are done.
So we assume at least one of the additional detectors is not in $N(v_2)$; due to symmetry, without loss of generality let $v_{24} \in S$.
Because $dom(v_3) = 3$ by hypothesis, $|\{v_{12},v_{13},v_{14}\} \cap S| = 1$.
If $v_{12} \in S \lor v_{13} \in S$ then $\{v_2,v_3\}$ are not distinguished, requiring $dom(v_2) = 3$, and we are done.
Similarly, if $v_{14} \in S$ then $\{v_3,v_4\}$ are not distinguished, requiring $dom(v_4) \ge 3$, and we are done.
In any case we see that $sh(x) \le \frac{10}{3}$.

\begin{wrapfigure}[10]{r}{0.3\textwidth}
    \vspace{-0.5em}
    \centering
    \includegraphics[width=0.25\textwidth]{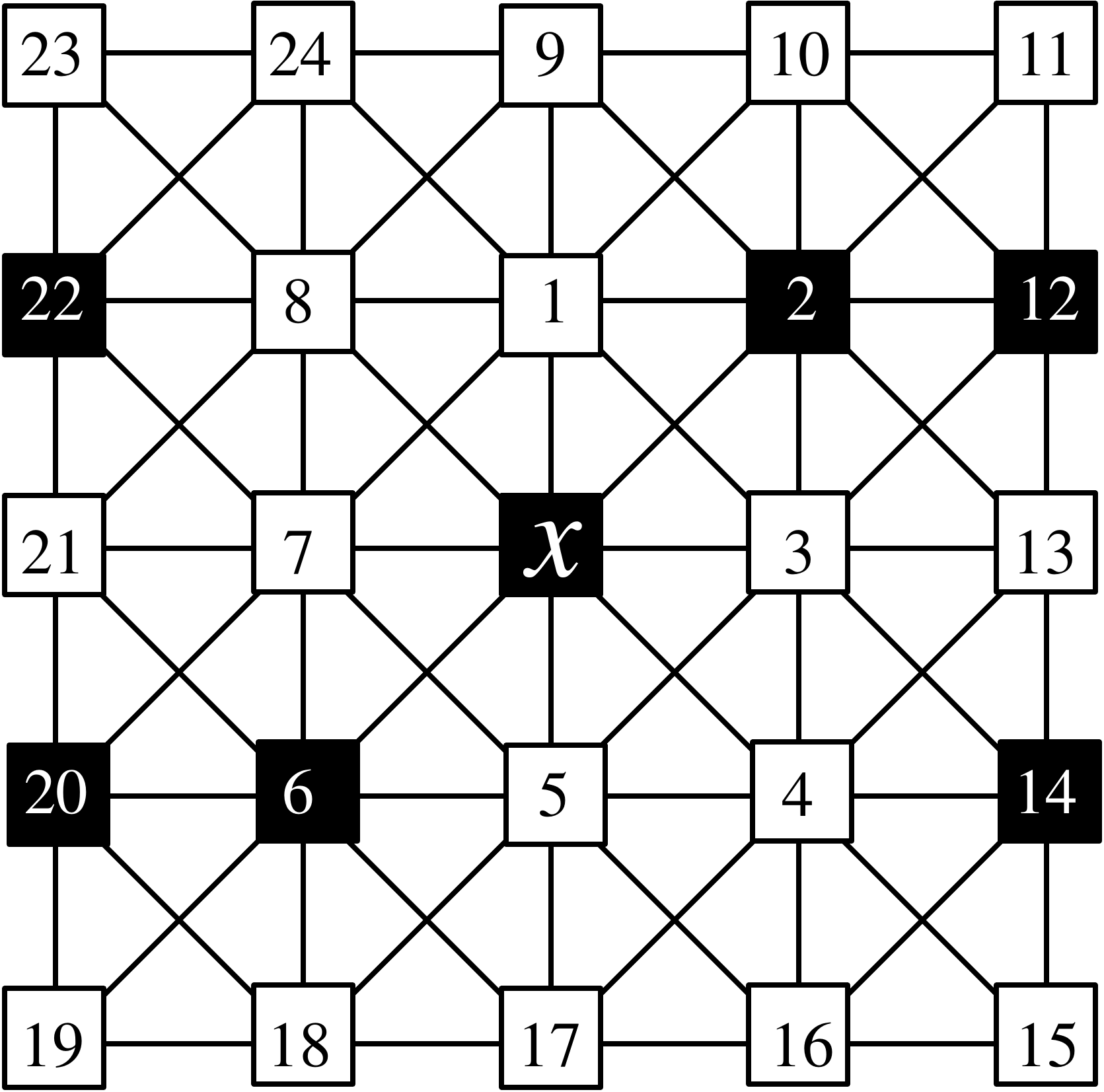}
    \caption{Final configuration of Case~5. $sh(x) > \frac{10}{3}$}
    \label{fig:red-old-case-5-big}
\end{wrapfigure}

If $dom(v_1) = 2 \land dom(v_3) < 4$ then $\{v_1, v_3\}$ are impossible to distinguish and thus the configuration is invalid.
Otherwise we have the final possibility: $dom(v_1) = 2 \le 4 = dom(v_3)$.
Again note that $\{v_5,v_7\}$ are not distinguished, so at least one of them must be at least 3-dominated.
Because $dom(v_1) = 2$, $\{v_9,v_{10},v_{24}\} \cap S = \varnothing$.
We now consider the 2 additional detectors in $N(v_3)$.
If $\{v_{12},v_{13}\} \subseteq S$ then $v_{14} \notin S$.
We observe that $\{v_2,v_3\}$ are not distinguished, requiring $v_{11} \in S$, making $dom(v_2) = 4$.
Thus, $sh(x) \le \frac{2}{4} + \frac{1}{3} + \frac{5}{2} = \frac{10}{3}$.
Otherwise if $\{v_{13},v_{14}\} \subseteq S$ then $v_{12} \notin S$.
If $v_{11} \in S$ then we have at least 4 vertices which are at least 4-dominated and so $sh(x) \le \frac{4}{3} + \frac{4}{2} = \frac{10}{3}$; so we consider only when $v_{11} \notin S$.
We see that $\{v_2,v_4\}$ are not distinguished, requiring $dom(v_4) \ge 4$.
Thus, $sh(x) \le \frac{2}{4} + \frac{1}{3} + \frac{5}{2} = \frac{10}{3}$.
Otherwise we have the final possible sub-configuration: $\{v_{12},v_{14}\} \subseteq S$, which means $v_{13} \notin S$.
We now look to vertex pair $(v_5,v_7)$, which is symmetric to the the pair $(v_1,v_3)$ we have thus far been using.
We applying the same arguments to $(v_5,v_7)$: consider $A = \{v_{16},v_{17},v_{18},v_{20},v_{21},v_{22}\} \cap S$.
By symmetry, any value of $A$ other than $\{v_{16},v_{18}\}$ or $\{v_{20},v_{22}\}$ results in $sh(x) \le \frac{10}{3}$.
If $A = \{v_{16},v_{18}\}$ then $dom(v_4) \ge 3$ and we are done, so we assume $A = \{v_{20},v_{22}\}$.
Thus, $dom(v_3) = dom(v_7) = 4$.
We see that if $\{v_{11},v_{15},v_{19},v_{23}\} \cap S \neq \varnothing$ then we have an additional neighbor of $x$ which is at least 3-dominated, so $sh(x) \le \frac{2}{4} + \frac{1}{3} + \frac{5}{2} = \frac{10}{3}$.
So we consider only the event when $\{v_{11},v_{15},v_{19},v_{23}\} \cap S = \varnothing$ and arrive at the sub-configuration in Figure~\ref{fig:red-old-case-5-big}.

Figure~\ref{fig:red-old-case-5-big} represents the final possible sub-configuration for Case~5, and we find that $sh(x) = \frac{7}{2} > \frac{10}{3}$.
We observe that $x$ has two neighbors, $\{v_2,v_6\}$, of ``bent" shape, which in Case~1 was found to have share no more than $\frac{13}{4}$.
Additionally, these two vertices are only adjacent to two other detector vertices $\{v_{12},v_{20}\} \subseteq S$, and are connected via a horizontal edge, which Cases~1--2 have proven to have share no more than $\frac{13}{4}$.
Thus $\{v_2,v_6\}$ are not adjacent to any other vertices with share exceeding $\frac{10}{3}$ and so are safe to use with averaging.
Consider the average share of $\{x,v_2,v_6\}$: $avg \le \frac{1}{3} \left[ \frac{7}{2} + \frac{13}{4} + \frac{13}{4} \right] = \frac{10}{3}$.

\FloatBarrier
\newcase{6}

\begin{wrapfigure}{r}{0.3\textwidth}
    \centering
    \includegraphics[width=0.25\textwidth]{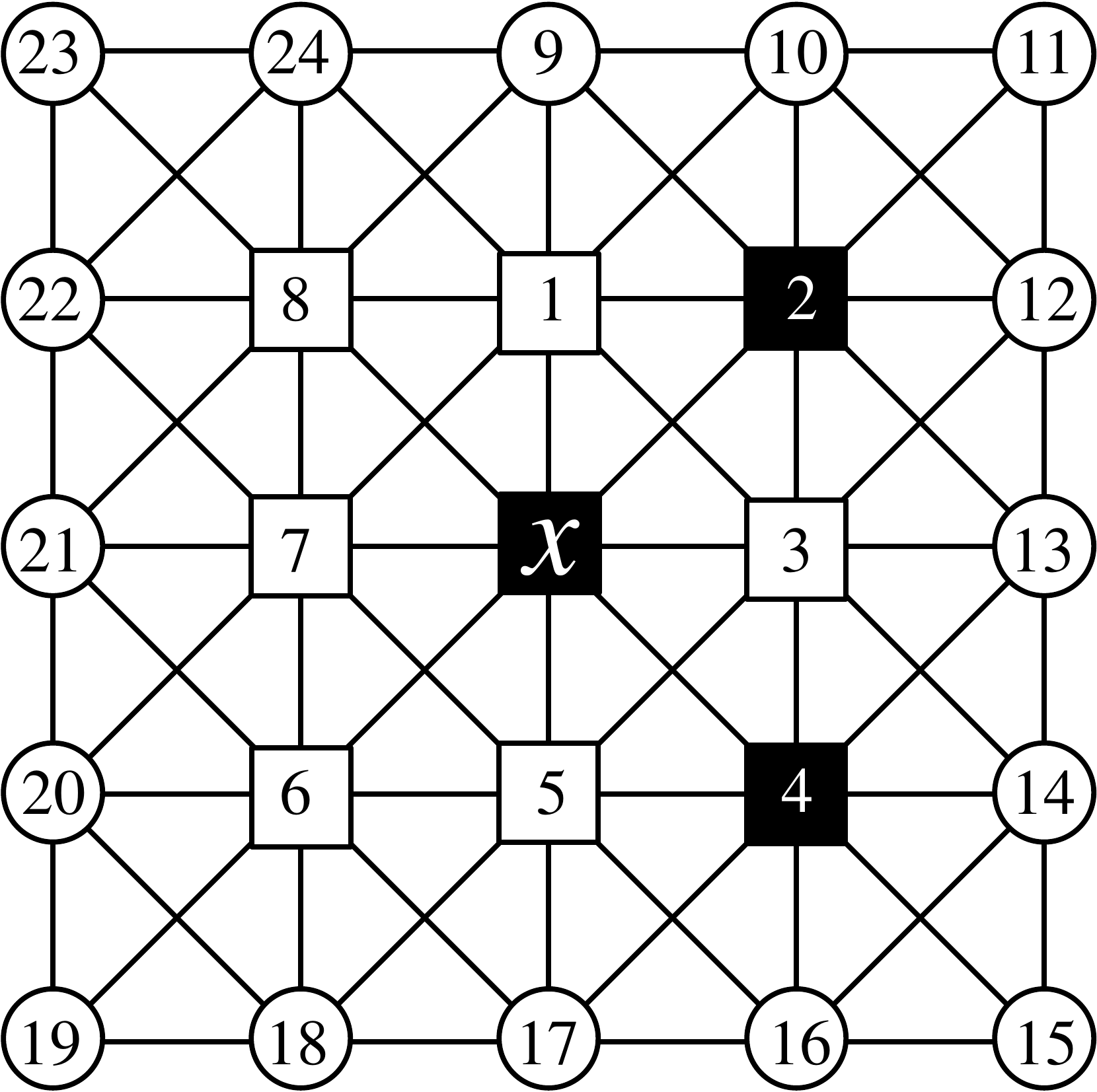}
    \caption{Case 6}
    \label{fig:red-old-case-6}
\end{wrapfigure}

In the final case, we consider when $x$ is dominated by two neighbors in a ``wedge" shape like Figure~\ref{fig:red-old-case-6}.
Applying the cumulative results of Cases~1--5 we see that if $\{v_1,v_3,v_5,v_6,v_7,v_8\} \cap S \neq \varnothing$ then the average share is no more than $\frac{10}{3}$.
We will proceed by expanding every possible state for the neighborhood of $v_7$.
Vertex $v_7$ must be at least 2-dominated, so we consider all sub-cases for $2 \leq dom(v_7) \le 4$.

First, we consider when $dom(v_7) = 4$, which requires $\{v_{20},v_{21},v_{22}\} \subseteq S$.
This implies that $\forall w \in \{v_6,v_7,v_8\},\; dom(w) \ge 3$.
Additionally, $dom(v_3) \ge 3$ by hypothesis.
Thus, $sh(x) \le \frac{4}{3} + \frac{4}{2} = \frac{10}{3}$.

There are two non-isomorphic ways to have $dom(v_7) = 3$: one with $\{v_{21}, v_{22}\} \subseteq S$ and the other with $\{v_{20}, v_{22}\} \subseteq S$.
In either way, $\{v_6,v_7\}$ are not distinguished, nor are $\{v_7,v_8\}$.
This requires $dom(v_6) \ge 3$ and $dom(v_8) \ge 3$.
Then $\forall w \in \{v_3,v_6,v_7,v_8\},\; dom(w) \ge 3$.
Thus, $sh(x) \le \frac{4}{3} + \frac{4}{2} = \frac{10}{3}$.

Finally, there are two non-isomorphic ways to have $dom(v_7) = 2$: one with $v_{22} \in S$ and the other with $v_{21} \in S$.  
If we have $dom(v_7) = 2$ with $v_{22} \in S$, we find that $\{v_7,v_8\}$ are not distinguished, requiring $dom(v_8) \ge 4$, which in turn requires $dom(v_1) \ge 3$.
Currently, we have three neighbors of $x$ which are at least 3-dominated; if we can find a fourth, then $sh(x) \le \frac{4}{3} + \frac{4}{2} = \frac{10}{3}$ and we will be done.
Thus, if $\exists w \in \{v_4, v_5, v_6\}$ such that $dom(w) \ge 3$ then we are done, so we assume $\forall w \in \{v_4, v_5, v_6\},\; dom(w) = 2$.
Vertex $v_5$ is already dominated by 2 neighbors, so $\{v_{16},v_{17},v_{18}\} \cap S = \varnothing$.
We see that $\{v_3,v_5\}$ are not distinguished, requiring $dom(v_3) \ge 4$.
Therefore we have that $dom(v_8) \ge 4$, $dom(v_3) \ge 4$, and $dom(v_1) \ge 3$, so $sh(x) \le \frac{2}{4} + \frac{1}{3} + \frac{5}{2} = \frac{10}{3}$.
If we have $dom(v_7) = 2$ with $v_{21} \in S$, then we observe that $\{v_6,v_7\}$ are not distinguished, nor are $\{v_7,v_8\}$.
This requires $dom(v_6) \ge 4$ and $dom(v_8) \ge 4$, which in turn makes $dom(v_1) \ge 3$ and $dom(v_5) \ge 3$.
As we already have that $dom(v_3) \ge 3$, $sh(x) \le \frac{2}{4} + \frac{3}{3} + \frac{3}{2} = 3 < \frac{10}{3}$.
\cendproof

From Theorem~\ref{theo:main-sh}, we have that the average share for any RED:OLD set on $K$ is no more than $\frac{10}{3}$.
Thus, $\frac{3}{10}$ is a lower bound for RED:OLD\%($K$).

\begin{corollary}
RED:OLD\%($K$) $\ge \frac{3}{10}$
\end{corollary}

\begin{figure}[ht]
    \centering
    \begin{tabular}{cc}
        \begin{tabular}{c}
           \includegraphics[width=0.34\textwidth]{fig/soln-1-both.pdf} \\ (a) \vspace{1.5em} \\
           \includegraphics[width=0.34\textwidth]{fig/soln-2-both.pdf} \\ (b)
        \end{tabular} &
        \begin{tabular}{c}
            \includegraphics[width=0.55\textwidth]{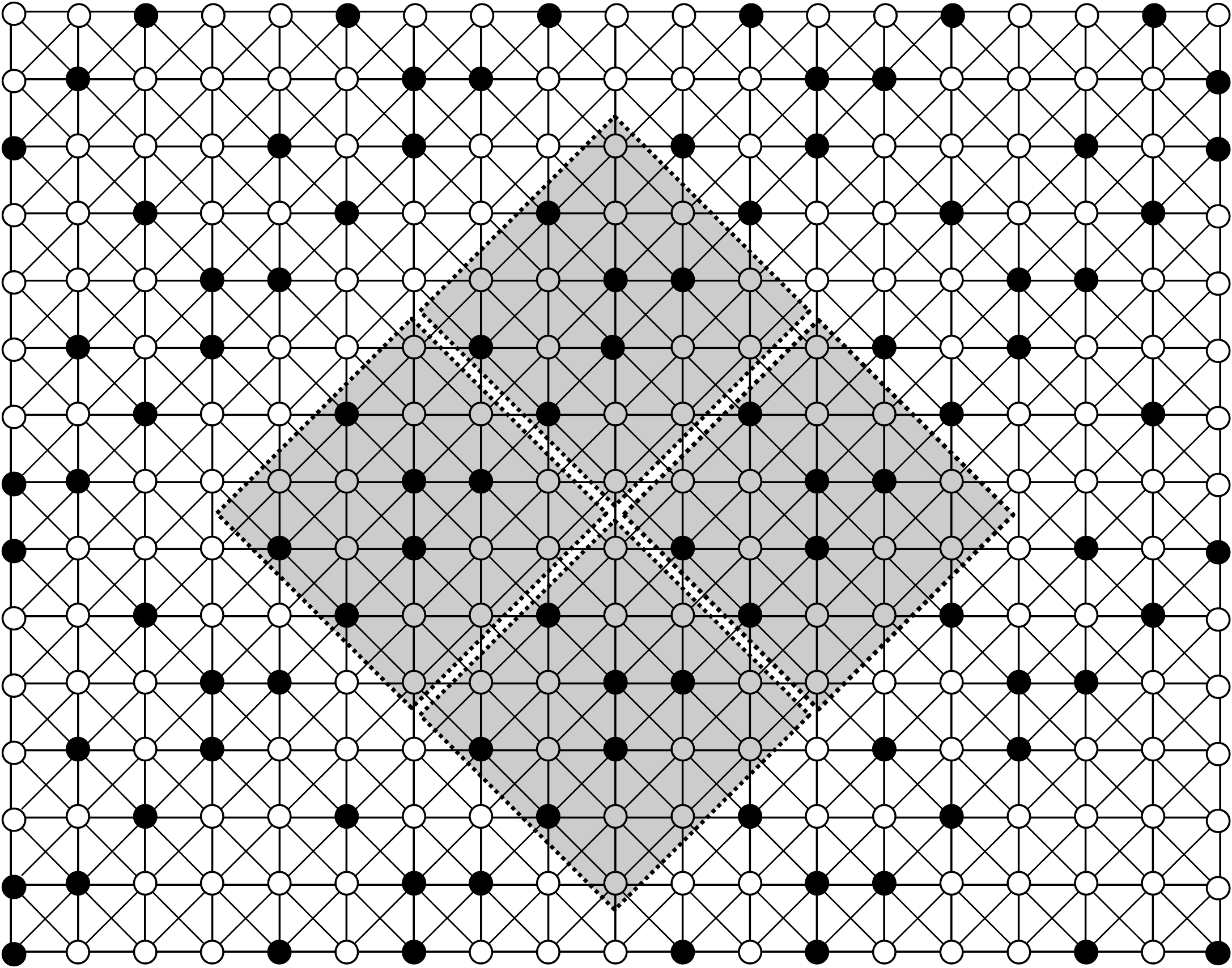} \\ (c)
        \end{tabular}
    \end{tabular}
    \caption{RED:ICs on $K$ with density $\frac{1}{3}$}
    \label{fig:red-ic-soln}
\end{figure}

\section{Bounds for $\textrm{RED:IC\%}(K)$}\label{sec:red-ic}

\subsection{Upper Bound}

We have found three RED:ICs on $K$, with density $\frac{1}{3}$, as shown in Figure~\ref{fig:red-ic-soln}.
Solutions (a) and (b) are identical to the RED:OLD solutions in Figure~\ref{fig:soln}---that is, these two solutions are both RED:OLD and RED:IC.
However, Figure~\ref{fig:red-ic-soln}~(c) is not a RED:OLD set due to some of the vertices being only 1-dominated if open neighborhoods were used.
One can easily verify that all of these solutions are RED:ICs via Theorem~\ref{theo:red-det-sys} with $k=1$ and $R_{det}(v) = N[v]$.  
Therefore, we have an upper bound for the optimal density: RED:IC\%($K$) $\le \frac{1}{3}$.

\FloatBarrier
\subsection{Lower Bound}

Let $S \subseteq V(G)$ be a RED:IC on $K$, and let the notations $sh$, $dom$, and $k$-dominated implicitly use $S$.
These will persist for the rest of Section~\thesubsection.

$S$ is a 1-redundant detection system, so Theorem~\ref{theo:red-det-sys} gives us that every vertex is at least 2-dominated.
Additionally, $K$ is 8-regular, so $|N[v]| = 9$.
Therefore, $\forall x \in S,\; sh(x) \le \frac{9}{2}$.
Thus, we have a simple lower bound: RED:IC\%($K$) $\ge \frac{2}{9}$.

\begin{theorem}\label{theo:main-sh-red-ic}
The share of any given detector in $S$ is at most $\frac{11}{3}$.
\end{theorem}

\cbeginproof
Just as in the previous section for proving a lower bound for RED:OLD\%($K$), we consider a detector $x \in S$ and its local region of $K$.
Vertex $x$ must be at least 2-dominated and dominates itself, so we consider the following two cases for the other detector being horizontal/vertical to $x$ or being diagonal to $x$, and show $sh(x) \le \frac{11}{3}$.

\FloatBarrier
\begin{wrapfigure}[4]{r}{0.2\textwidth}
    \vspace{-1.8em}
    \centering
    \includegraphics[width=0.15\textwidth]{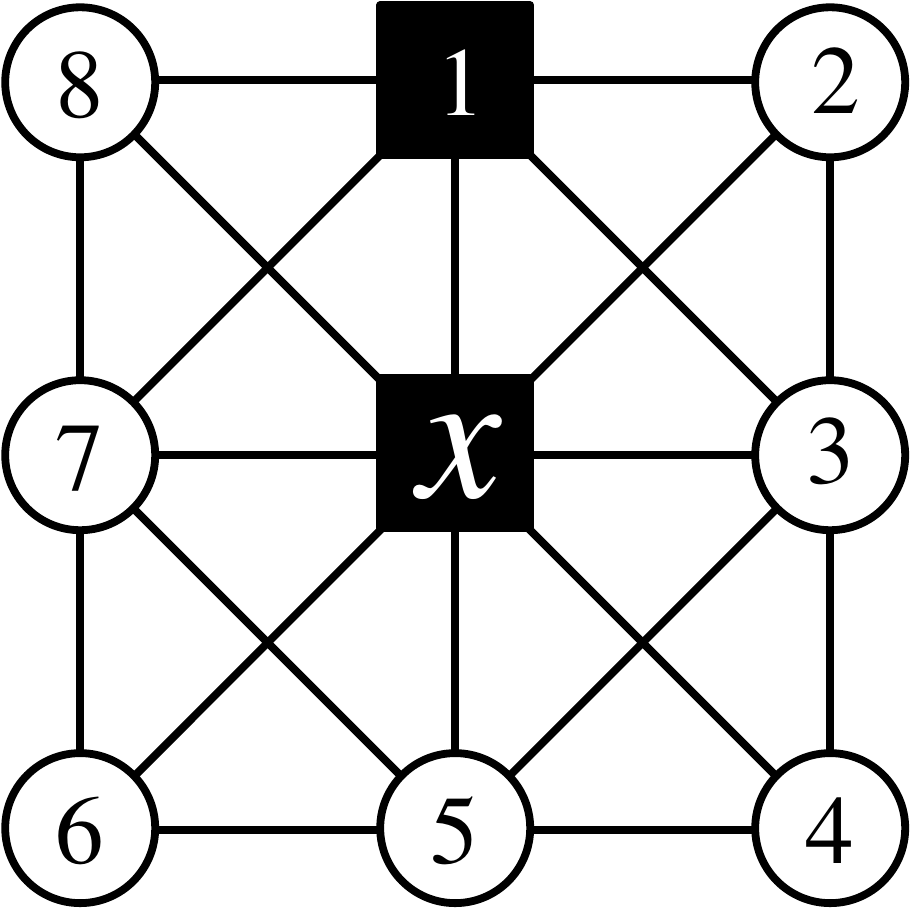}
    \caption{Case~1}
    \label{fig:red-ic-case-1}
\end{wrapfigure}
\newcase{1}

First, we consider when $v_1 \in S$ and everything else is unknown, as in Figure~\ref{fig:red-ic-case-1}.
Consider $D = \{x,v_1\}$ and $A = \{x,v_1,v_2,v_3,v_7,v_8\}$.
By Lemma~\ref{lem:common}, $sh[A](x) \le \frac{6}{2+1} = 2$.
We can assume the other three vertices $\{v_4,v_5,v_6\}$ are at least 2-dominated.
Therefore, $sh(x) \le 2 + \frac{3}{2} = \frac{7}{2} < \frac{11}{3}$.

\FloatBarrier
\newcase{2}

\begin{wrapfigure}[8]{r}{0.2\textwidth}
    \vspace{-0.5em}
    \centering
    \includegraphics[width=0.15\textwidth]{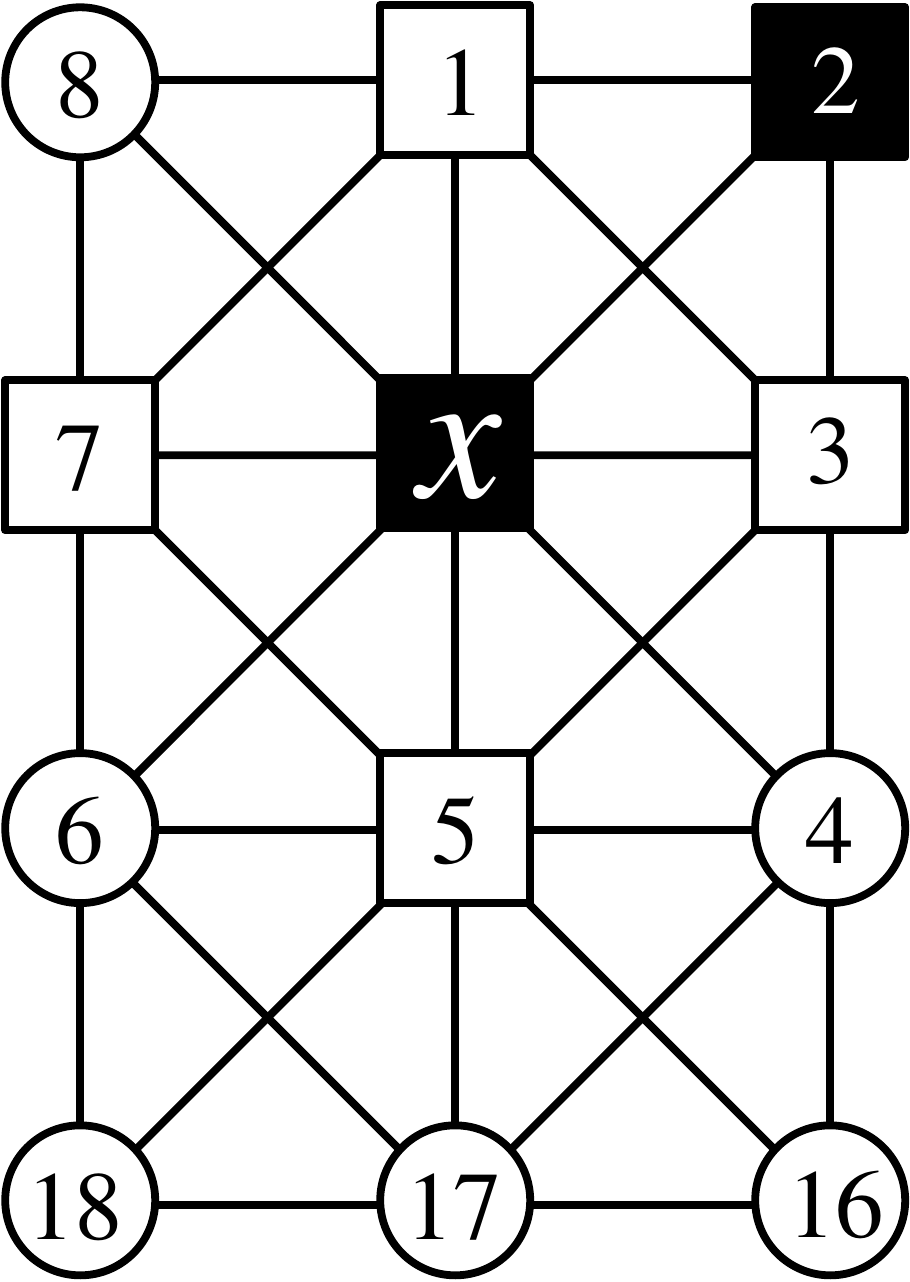}
    \caption{Case~2}
    \label{fig:red-ic-case-2}
\end{wrapfigure}
Now we consider the only other case, when $v_2 \in S$.
From the results of Case~1, we know that if $\{v_1,v_3,v_5,v_7\} \cap S \neq \varnothing$ then $sh(x) \le \frac{11}{3}$; so we assume $\{v_1,v_3,v_5,v_7\} \cap S = \varnothing$, as shown in Figure~\ref{fig:red-ic-case-2}.
Let $D = \{x,v_2\}$ and $A = \{x,v_1,v_2,v_3\}$; by Lemma~\ref{lem:common}, $sh[A](x) \le \frac{4}{2+1} = \frac{4}{3}$.
Suppose $\forall w \in \{v_4,v_5,v_6\},\; dom(w) = 2$.

We will focus on $v_5$, which by hypothesis is 2-dominated, so there must be a detector in $N[v_5]$ other than $x$.
However, if the second detector is in $N[v_4] \cap N[v_5]$ then $\{v_4,v_5\}$ will share them both and therefore not be distinguished; thus, $\{v_4,v_{16},v_{17}\} \cap S = \varnothing$.
Similarly, if the second detector is in $N[v_5] \cap N[v_6]$ then $\{v_5,v_6\}$ will not be distinguished; thus, $\{v_6,v_{17},v_{18}\} \cap S = \varnothing$.
Thus, we have $\{v_4,v_6,v_{16},v_{17},v_{18}\} \cap S = \varnothing$, which results in $dom(v_5) =1$, a contradiction.
Therefore, $\exists w \in \{v_4,v_5,v_6\}$ such that $dom(w) \ge 3$, making $sh(x) \le \frac{4}{3} + \frac{1}{3} + \frac{4}{2} = \frac{11}{3}$ and we are done.
\cendproof

From Theorem~\ref{theo:main-sh-red-ic}, we have that the average share for any RED:IC on $K$ is no more than $\frac{11}{3}$.
Thus, $\frac{3}{11}$ is a lower bound for RED:IC\%($K$).

\begin{corollary}
RED:IC\%($K$) $\ge \frac{3}{11}$
\end{corollary}

\FloatBarrier
\section{Conclusion}

We have now demonstrated that $\frac{3}{10} \le$ RED:OLD\%($K$) $\le \frac{1}{3}$ and $\frac{3}{11} \le$ RED:IC\%($K$) $\le \frac{1}{3}$.
Interestingly, two of the solutions achieving the upper bound of $\frac{1}{3}$ for RED:OLD and RED:IC were the same subset: (a) and (b) from Figures \ref{fig:soln} and \ref{fig:red-ic-soln}.
For many graphs, RED:IC\%($G$) $<$ RED:OLD\%($G$), but in this case, the upper bounds we obtained are equal.
Although it is possible that RED:IC\%($K$) $=$ RED:OLD\%($K$), it would not be surprising if the provided RED:ICs were found to be sub-optimal.

% \section*{Acknowledgement} 

\bibliographystyle{acm}
\bibliography{refs}

\begin{thebibliography}{10}

\bibitem{watchsys-1}
{\sc Auger, D., Charon, I., Hudry, O., and Lobstein, A.}
\newblock Watching systems in graphs: an extension of identifying codes.
\newblock {\em Discrete Applied Mathematics 161}, 12 (2013), 1674--1685.

\bibitem{watchsys-2}
{\sc Auger, D., Charon, I., Hudry, O., and Lobstein, A.}
\newblock Maximum size of a minimum watching system and the graphs achieving
  the bound.
\newblock {\em Discrete Applied Mathematics 164\/} (2014), 20--33.

\bibitem{chhl}
{\sc Charon, I., Honkala, I., Hudry, O., and Lobstein, A.}
\newblock The minimum density of an identifying code in the king lattice.
\newblock {\em Discrete Mathematics 276}, 1--3 (2004), 95--109.

\bibitem{chl}
{\sc Charon, I., Hudry, O., and Lobstein, A.}
\newblock Identifying codes with small radius in some infinite regular graphs.
\newblock {\em The Electronic Journal of Combinatorics 9}, 1 (2002), R11.

\bibitem{chlz}
{\sc Cohen, G.~D., Honkala, L., and Lobstein, A.}
\newblock On codes identifying vertices in the two-dimensional square lattice
  with diagonals.
\newblock {\em IEEE Transactions on Computers 50}, 2 (2001), 174--176.

\bibitem{dantas}
{\sc Dantas, R., Havet, F., and Sampaio, R.~M.}
\newblock Minimum density of identifying codes of king grids.
\newblock {\em Discrete Mathematics 341}, 10 (2018), 2708--2719.

\bibitem{flp}
{\sc Foucaud, F., Laihonen, T., and Parreau, A.}
\newblock An improved lower bound for $(1,\leq 2)$-identifying codes in the
  king grid.
\newblock {\em Advances in Mathematics of Communications 8}, 1 (2014), 35--52.

\bibitem{hl2}
{\sc Honkala, I., and Laihonen, T.}
\newblock On locating--dominating sets in infinite grids.
\newblock {\em European Journal of Combinatorics 27}, 2 (2006), 218--227.

\bibitem{hl3}
{\sc Honkala, I., and Laihonen, T.}
\newblock On identifying codes in the king grid that are robust against edge
  deletions.
\newblock {\em The Electronic Journal of Combinatorics\/} (2008), R3--R3.

\bibitem{our13}
{\sc Jean, D.~C., and Seo, S.~J.}
\newblock Fault-tolerant identifying codes in special classes of graphs.
\newblock {\em Discussiones Mathematicae: Graph Theory (in press)\/} (2021).

\bibitem{ourtri}
{\sc Jean, D.~C., and Seo, S.~J.}
\newblock {Optimal Error-Detecting Open-Locating-Dominating Set on the Infinite
  Triangular Grid}.
\newblock {\em Discussiones Mathematicae: Graph Theory (in press)\/} (2021).

\bibitem{karpovsky}
{\sc Karpovsky, M.~G., Chakrabarty, K., and Levitin, L.~B.}
\newblock {On a new class of codes for identifying vertices in graphs}.
\newblock {\em IEEE Transactions on Information Theory 44}, 2 (1998), 599--611.

\bibitem{oldtri}
{\sc Kincaid, R., Oldham, A., and Yu, G.}
\newblock Optimal open-locating-dominating sets in infinite triangular grids.
\newblock {\em Discrete Applied Mathematics 193\/} (2015), 139--144.

\bibitem{dombib}
{\sc lobstein, A.}
\newblock {Watching Systems, Identifying, Locating-Dominating and
  Discriminating Codes in Graphs}.
\newblock \url{https://www.lri.fr/\%7elobstein/debutBIBidetlocdom.pdf}.

\bibitem{p2}
{\sc Pelto, M.}
\newblock On locating-dominating codes for locating large numbers of vertices
  in the infinite king grid.
\newblock {\em Australasian Journal of Combinatorics 50\/} (2011), 127--139.

\bibitem{p5}
{\sc Pelto, M.}
\newblock Optimal $(r,\le3)$-locating–dominating codes in the infinite king
  grid.
\newblock {\em Discrete Applied Mathematics 161}, 16 (2013), 2597--2603.

\bibitem{p6}
{\sc Pelto, M.}
\newblock Optimal identifying codes in the infinite 3-dimensional king grid.
\newblock {\em European Journal of Combinatorics 36\/} (feb 2014), 641–659.

\bibitem{p4}
{\sc Pelto, M., et~al.}
\newblock On identifying and locating-dominating codes in the infinite king
  grid.
\newblock {\em TUCS Dissertations No. 155\/} (2012).

\bibitem{oldking}
{\sc Seo, S.}
\newblock {Open-locating-dominating sets in the infinite king grid}.
\newblock {\em Journal of Combinatorial Mathematics and Combinatorial Computing
  104\/} (2018), 31--47.

\bibitem{ft-old-cubic}
{\sc Seo, S.~J.}
\newblock {Fault-tolerant detectors for distinguishing sets in cubic graphs}.
\newblock {\em Discrete Applied Mathematics 293\/} (2021), 25--33.

\bibitem{old}
{\sc Seo, S.~J., and Slater, P.~J.}
\newblock {Open neighborhood locating dominating sets.}
\newblock {\em Australasian Journal of Combinatorics 46\/} (2010), 109--120.

\bibitem{tmb}
{\sc Seo, S.~J., and Slater, P.~J.}
\newblock Graphical parameters for classes of tumbling block graphs.
\newblock {\em Congressus Numerantium 213\/} (2012), 155--168.

\bibitem{ftsets}
{\sc Seo, S.~J., and Slater, P.~J.}
\newblock {Fault tolerant detectors for distinguishing sets in graphs}.
\newblock {\em Discussiones Mathematicae: Graph Theory 35\/} (2015), 797--818.

\bibitem{dom-loc-acyclic}
{\sc Slater, P.~J.}
\newblock {Domination and location in acyclic graphs}.
\newblock {\em Networks 17}, 1 (1987), 55--64.

\bibitem{ftld}
{\sc Slater, P.~J.}
\newblock {Fault-tolerant locating-dominating sets}.
\newblock {\em Discrete Mathematics 249}, 1--3 (2002), 179--189.

\end{thebibliography}

\end{document}